\documentclass[11pt,a4paper]{article}

\usepackage[T1]{fontenc}
\usepackage[utf8]{inputenc}
\usepackage{lmodern}
\usepackage{microtype}
\usepackage{mathtools,amsthm}
\usepackage{amssymb,amsfonts}
\usepackage{mathrsfs}
\usepackage{enumitem}
\usepackage{xcolor}
\usepackage{tikz}
\usetikzlibrary{arrows.meta}
\usepackage{placeins}
\usepackage{hyperref}
\usepackage[nameinlink,capitalize,noabbrev]{cleveref}

\hypersetup{
  colorlinks=true,
  linkcolor=blue!60!black,
  citecolor=green!40!black,
  urlcolor=blue!60!black
}

\allowdisplaybreaks[2]
\numberwithin{equation}{section}

\newtheorem{theorem}{Theorem}[section]
\newtheorem{proposition}[theorem]{Proposition}
\newtheorem{lemma}[theorem]{Lemma}
\newtheorem{corollary}[theorem]{Corollary}
\newtheorem{auxiliarylemma}[theorem]{Lemma}
\theoremstyle{definition}

\newtheorem{introdefinition}{Definition}[section]
\crefname{introdefinition}{definition}{definitions}
\Crefname{introdefinition}{Definition}{Definitions}

\theoremstyle{remark}
\newtheorem{remark}[theorem]{Remark}

\newcommand{\T}{\mathbb T}
\newcommand{\R}{\mathbb R}

\newcommand{\N}{\mathbb N}
\newcommand{\Z}{\mathbb Z}
\newcommand{\cK}{\mathcal K}
\newcommand{\cA}{\mathcal A}
\newcommand{\cS}{\mathcal S}

\newcommand{\ip}[2]{\langle #1,#2\rangle}
\newcommand{\norm}[1]{\left\lVert #1\right\rVert}
\newcommand{\abs}[1]{\left\lvert #1\right\rvert}
\newcommand{\tr}{\operatorname{tr}}
\newcommand{\dist}{\operatorname{dist}}
\newcommand{\Graph}{\operatorname{Graph}}

\newcommand{\op}{\mathrm{op}}
\newcommand{\HS}{\mathfrak I_2}
\newcommand{\Ione}{\mathfrak I_1}

\newcommand{\dd}{\,\mathrm d}

\usepackage[textwidth=145mm,textheight=225mm,centering]{geometry}
\title{\Large\bfseries Global attractors for the weakly damped forced periodic KdV equation in $H^{-1}$}
\hypersetup{pdfauthor={Minjie Shan}}
\date{}
\author{\footnotesize\MakeUppercase{Minjie Shan}}

\begin{document}
\maketitle
\begin{abstract}
The long-time dynamics of the weakly damped periodic Korteweg--de Vries equation are investigated for mean-zero initial data in $H^{-1}$ and time-independent forcing in $L^2$. Primitive trajectories are obtained as strong $H^{-1}$ limits of global smooth solutions that evolve without restarts, with convergence locally uniform in time. Fourier truncations of the initial data and forcing provide such trajectories from every initial state. A strong trajectory attractor is established for this class. To obtain the strict semiflow property in phase space without assuming uniqueness, we enlarge the trajectory class by allowing piecewise restarts of the smooth approximants. The total $H^{-1}$ jump size is required to tend to zero on every finite time interval. We prove that  the resulting defect trajectories generate a strict multivalued semiflow admitting a compact global attractor in the strong $H^{-1}$ topology. This attractor is further shown to be bounded in $L^2$ and to coincide with the classical $L^2$ attractor. The proof combines dissipative estimates for the perturbation determinant associated with the KdV Lax operator and the Birkhoff coordinate equations. The determinant estimates provide spatial compactness, while the coordinate equations yield equicontinuity in time on compact subsets of the phase space.
\end{abstract}
\begin{center}
\begin{minipage}{0.94\textwidth}
\small
\textbf{Keywords:} KdV equation; weak damping; global attractor;
perturbation determinant; Birkhoff coordinates.\par\smallskip
\textbf{MSC 2020:} 35Q53, 35B41, 37L30, 47H20.
\end{minipage}
\end{center}
\medskip

\section{Introduction}\label{sec:introduction}

We consider the weakly damped and forced periodic Korteweg--de Vries equation
\begin{equation}\label{eq:fdkdv}
 \begin{cases}
 u_t+u_{xxx}-3\partial_x(u^2)+\gamma u=f,
     &(t,x)\in\mathbb R\times\mathbb T,\\
 u(0,x)=u_0(x)\in H^{-1}_0(\mathbb T),
 \end{cases}
\end{equation}
where $\mathbb T=\mathbb R/(2\pi\mathbb Z)$, $\gamma>0$, and
$f\in L^2_0(\mathbb T)$ is independent of time. All states are real-valued,
and
\[
 H^s_0(\mathbb T)
 =\{q\in H^s(\mathbb T;\mathbb R):
       \langle q,1\rangle_{\mathcal D',\mathcal D}=0\},
 \qquad s\in\mathbb R.
\]
The distributional pairing with $1$ defines the zero-mean condition
at every Sobolev regularity. For a smooth solution, we write 
\[
 m(t)=\frac1{2\pi}\int_{\mathbb T}u(t,x)\,dx.
\]
Then, integration of
\eqref{eq:fdkdv} over $\mathbb T$ gives
\[
 m'(t)+\gamma m(t)=0,  \qquad
\text{and}  \qquad
 m(t)=e^{-\gamma t}m(0).
\]
Thus the evolution preserves zero spatial mean. The same condition
passes to the limits of smooth solutions considered below.

The analysis of \eqref{eq:fdkdv} in $H^{-1}$ presents substantial difficulties. At this regularity, $u$ need not be a function, and the quadratic nonlinearity $\partial_x(u^2)$ is not well defined in general. Moreover, even strong convergence in $H^{-1}$ does not suffice to justify passage to the limit in this term. Indeed, setting $r_N(x)=\cos(Nx)+\cos((N+1)x)$, we have
\begin{equation}\label{eq:quadratic-discontinuity}
 r_N\longrightarrow0\quad\text{in }H^{-1},
 \qquad
 \partial_x(r_N^2)\longrightarrow-\sin x
 \quad\text{in }\mathcal D'(\mathbb T).
\end{equation}
Nevertheless, the integrable structure allows the conservative KdV flow
to be extended continuously to $H^{-1}$
\cite{KappelerTopalov2006,KillipVisan2019}. In Birkhoff coordinates,
the evolution consists of rotations whose frequencies are determined
by the conserved actions. Under the perturbation $f-\gamma u$,
however, these actions are generally no longer conserved. Moreover,
continuity of the conservative flow does not by itself provide an estimate for the
difference between two forced and damped solutions. It therefore does
not imply that smooth approximations of the same rough initial datum
have a common limit.

A forward \emph{trajectory} is a continuous path
$u:[0,\infty)\to H^{-1}_0(\mathbb T)$ with $u(0)=u_0$.
Since the initial state may be a distribution, we must specify which
paths are admissible for \eqref{eq:fdkdv}. We use limits of smooth
solutions with Fourier-truncated forcing and smooth initial states
converging in $H^{-1}$. In particular, Fourier truncations of the initial
datum provide the approximations used to construct trajectories.
The resulting primitive and defect classes are introduced in
\Cref{def:primitive,def:defect-trajectory}. This framework describes
long-time dynamics both in the space of forward paths and in the
phase space $H^{-1}_0$. We also investigate the spatial regularity of
bounded trajectories defined for all $t\in\mathbb R$.

\subsection{Notation}

 We use normalized measure $d\mu=(2\pi)^{-1}dx$ on $\mathbb T$ and set
\[
 \widehat q(k)=\frac{1}{2\pi}\langle q,e^{-ikx}\rangle,\qquad
 \norm q_{H^s}^2=\sum_{k\in\mathbb Z}(1+k^2)^{s}
 |\widehat q(k)|^2.
\]
Here $\langle q,\phi\rangle$ is the distributional pairing extending
integration against $\phi$ with respect to $dx$.
For $\kappa\ge1$, we use the weighted $H^{-1}$ norm defined by
\begin{equation}\label{eq:weighted-Hminus1}
 \norm{q}_{H^{-1}_\kappa}^2
 =\sum_{k\in\mathbb Z}\frac{\abs{\widehat q(k)}^2}{k^2+4\kappa^2}.
\end{equation}
For each fixed $\kappa$, this norm is equivalent to the usual $H^{-1}$
norm. Let $P_{\le N}$ and $P_{>N}$ denote the Littlewood--Paley
projections onto $|k|\le N$ and $|k|>N$, respectively. For integers
$\nu\ge0$, set $f_\nu=P_{\le\nu}f$. All smooth approximants are
real-valued and have zero spatial mean. For $z=(z_1, z_2, \cdots, z_n, \cdots)$, define
\[
 \norm{z}_{h^s}^2=\sum_{n\ge1}n^{2s}|z_n|^2.
\]
The symbols $\Ione$ and $\HS$ denote the trace and Hilbert--Schmidt
classes on complexified $L^2(\mathbb T)$. Positive constants $c,C$ may
change from one estimate to the next; any relevant parameter dependence
is indicated. All trajectory convergence is locally uniform in time
with values in $H^{-1}$ unless another topology is specified.

\subsection{Main results}\label{1.2}

We now state the main results of this paper. For $\kappa\ge1$, let
\[
 R_0(\kappa)=(-\partial_x^2+\kappa^2)^{-1},\qquad
 A_\kappa(q)=R_0(\kappa)^{1/2}qR_0(\kappa)^{1/2}.
\]
For smooth $q$ the middle factor in $A_\kappa(q)$ is multiplication by $q$.
The operator extends by density to a Hilbert--Schmidt operator for
$q\in H^{-1}$. In the region $\norm{A_\kappa(q)}_{\op}<1$, we define
\begin{align}
 \alpha_\kappa(q)
 :=-\log\det\nolimits_{2}\bigl(I+A_\kappa(q)\bigr)
  =\sum_{\ell=2}^{\infty}\frac{(-1)^\ell}{\ell}
   \operatorname{tr}\bigl(A_\kappa(q)^\ell\bigr).
 \label{trace2aG}
\end{align}
Here $\det\nolimits_2$ is the regularized Fredholm determinant and the logarithm
is real.  

The following theorem provides the dissipative estimate needed for the construction of trajectories and attractors. This estimate is uniform in the Fourier cutoff parameter $\nu$ for the forcing and can be passed to both primitive and defect trajectories. The resulting absorbing bound and high-frequency tail estimates provide the spatial control needed to establish strong compactness of trajectories. Additionally, when applied to bounded complete trajectories, the estimate yields a uniform $L^2$ bound. This bound will be used in the subsequent identification of these trajectories with the classical $L^2$ dynamics.
\begin{theorem} \label{thm:intro-determinant}
Fix $\gamma>0$ and $f\in L^2_0$. Let $u$ be a global smooth solution of \eqref{eq:fdkdv} with force $f_\nu$ for some $\nu\ge0$. For each $R>0$ there is a threshold
$\kappa_R=\kappa(R,\gamma,\norm f_{H^{-1}})\ge1$ with the following
property.  If $\norm{u(0)}_{H^{-1}}\le R$,
then for all $\kappa\ge\kappa_R$ the determinant is defined along $u$
and satisfies
\begin{equation}\label{eq:intro-spectral-comparison}
 c\norm{u(t)}_{H^{-1}_\kappa}^2
 \le \kappa\alpha_\kappa(u(t))
 \le C\norm{u(t)}_{H^{-1}_\kappa}^2,\qquad t\ge0.
\end{equation}
Moreover, for every $0\le s\le t$,
\begin{equation}\label{eq:intro-dissipation}
 \kappa\alpha_\kappa(u(t))
 \le e^{-c\gamma(t-s)}\kappa\alpha_\kappa(u(s))
 +C\gamma^{-2}\norm{f}_{H^{-1}_\kappa}^2\bigl(1-e^{-c\gamma(t-s)}\bigr).
\end{equation}
The constants $c,C>0$ in both estimates are absolute.
In particular, the threshold $\kappa_R$  and estimates are uniform over the smooth
force cutoffs.
\end{theorem}

The proof combines \Cref{thm:forced-determinant}\textup{(iii)} with
\eqref{eq:det-weighted-norm}.

The dissipative estimate above is formulated for smooth solutions.
We next introduce trajectories in $H^{-1}_0$ through smooth approximation to study the long-time dynamics in path space.
We begin with smooth solutions that evolve from their initial states without restarts; their limiting paths are called primitive trajectories.

\begin{introdefinition} \label{def:primitive}
A path $u\in C_{\mathrm{loc}}([0,\infty);H^{-1}_0)$ is a
\emph{primitive trajectory} if there exist integers $\nu_j\to\infty$
and global smooth solutions $u_j$ of \eqref{eq:fdkdv} with forcing
$f_{\nu_j}$ and smooth initial states $q_j\to u(0)$ in $H^{-1}$ such that
\begin{equation}\label{eq:primitive-convergence}
 \lim_{j\to\infty}\sup_{0\le t\le T}
 \norm{u_j(t)-u(t)}_{H^{-1}}=0
 \qquad\text{for every }T>0.
\end{equation}
The same approximating sequence is used on every finite time interval,
and each $u_j$ evolves without restarts on $[0,\infty)$.
The class of all primitive trajectories is denoted by $\cK_0^+$.
\end{introdefinition}

To describe evolution in this trajectory class, we use time translation.
For $s\ge0$, define
\[
 (\mathcal{T}_s u)(t)=u(t+s),\qquad t\ge0.
\]
These operators satisfy $\mathcal{T}_0u=u$ and
$\mathcal{T}_s\mathcal{T}_r u=\mathcal{T}_{s+r}u$ for all $s,r\ge0$.
By autonomy of \eqref{eq:fdkdv}, $\mathcal{T}_s u_j$ is again a global
smooth solution with the same forcing as $u_j$.
Moreover, the convergence in \eqref{eq:primitive-convergence} is
preserved under time translation.
Thus $\mathcal{T}_s u\in\cK_0^+$ whenever $u\in\cK_0^+$.

The trajectory-attractor approach was developed in
\cite{ChepyzhovVishik1997,ChepyzhovVishik2002}.
Related frameworks for evolution equations without uniqueness
are provided by generalized semiflows \cite{Ball1997}
and evolutionary systems \cite{Cheskidov2009}.

We can now formulate attraction in path space. Here attraction is measured on each fixed time interval $[0,T]$ by comparing the translated trajectory $\mathcal{T}_s u$ with a single trajectory in the attracting set.
The following definition makes this precise. Here we use the strong trajectory topology as in
\cite[Chapter~II, Definition~4.1]{VishikChepyzhov2011} and
\cite[Definition~2.8 and Theorem~4.3]{ChepyzhovVishikZelik2011}.

\begin{introdefinition}
\label{def:trajectory-attractor}
Fix a class of forward trajectories closed under $\mathcal{T}_s$
for every $s\ge0$.
A strong trajectory attractor is a nonempty set $\cA$ of trajectories
from this class, compact in $C_{\mathrm{loc}}([0,\infty);H^{-1}_0)$
and strictly invariant under time translation:
\[
 \mathcal{T}_s\cA=\{\mathcal{T}_s u:u\in\cA\}=\cA,
 \qquad s\ge0.
\]
Writing $\mathcal K^+$ for the fixed trajectory class, we also require
\begin{equation}\label{eq:trajectory-attraction}
 \lim_{s\to\infty}
 \sup_{\substack{u\in\mathcal K^+\\u(0)\in B}}
 \inf_{v\in\cA}\norm{\mathcal T_su-v}_{C([0,T];H^{-1})}=0
\end{equation}
for every $T>0$ and every bounded $B\subset H^{-1}_0$.
\end{introdefinition}

 The next theorem establishes the existence of a primitive trajectory attractor and further describes its structure. To state this result, we first introduce the notion of a bounded complete primitive trajectory.

A path $U\in C(\mathbb R;H^{-1}_0)$ is called a
\emph{complete primitive trajectory} if
\[
 U(s+\cdot)|_{[0,\infty)}\in\cK_0^+
 \qquad\text{for every }s\in\mathbb R.
\]
Thus, the forward path starting at any real time is a primitive trajectory.
Such a trajectory is \emph{bounded} if
\[
 \sup_{t\in\mathbb R}\norm{U(t)}_{H^{-1}}<\infty.
\]

\begin{theorem}\label{thm:intro-primitive}
Let $\gamma>0$ and $f\in L^2_0(\mathbb T)$.
For every $u_0\in H^{-1}_0$, there exists a primitive trajectory
$u\in\cK_0^+$ with $u(0)=u_0$.
The class $\cK_0^+$ admits a unique strong trajectory attractor $\cA$
in the sense of \Cref{def:trajectory-attractor}.
Moreover,
\[
 \cA=\{U|_{[0,\infty)}:
 U\text{ is a bounded complete primitive trajectory}\}.
\]
If $f=0$, then $\cA$ consists of the zero trajectory.
\end{theorem}

The proof is given in Section~\ref{sec:primitive}.
Compactness of smooth solutions on finite time intervals yields primitive trajectories.
The trajectory attractor is then obtained from the compactness of large time translates of these trajectories.

The preceding result describes attraction in path space.
We now turn to the dynamics in the phase space $H^{-1}_0$.
Time translation follows the future of a fixed trajectory.
To obtain the strict semiflow property, we also need a trajectory class
that allows the first part of one trajectory to be joined to any
admissible continuation from the same state.

\begin{figure}[htbp]
  \centering
  \begin{tikzpicture}[x=1cm,y=1cm,>=Stealth,
      every node/.style={font=\small}]
    \draw[->,gray] (.35,.3)--(12.6,.3) node[right] {$t$};
    \draw[->,gray] (.35,.3)--(.35,4.15);
    \node[rotate=90,font=\footnotesize,text=gray] at (-.02,2.25)
      {state};
    \node[below] at (.65,.3) {$0$};
    \draw[gray!65,densely dotted] (5.5,.3)--(5.5,2.15);
    \node[below] at (5.5,.3) {$s$};

    \draw[blue!35!black,line width=.7pt]
      (.65,3.55) .. controls (2.1,3.75) and (3.6,2.65) .. (5.5,2.15)
      .. controls (7.4,1.65) and (9.65,1.0) .. (12.0,1.05);
    \draw[blue!75!black,line width=1.55pt,
      -{Stealth[length=2.3mm,width=1.7mm]}]
      (.65,3.55) .. controls (2.1,3.75) and (3.6,2.65) .. (5.5,2.15);

    \draw[teal!75!black,line width=1.55pt,
      -{Stealth[length=2.3mm,width=1.7mm]}]
      (5.5,2.15) .. controls (7.3,2.7) and (9.65,3.85) .. (12.0,3.55);
    \fill[blue!75!black] (5.5,2.15) circle[radius=2.6pt];
    \node[fill=white,inner sep=2pt] at (6.0,1.45) {$u(s)=v(0)$};
    \node[above,text=blue!75!black] at (11.55,1.1) {$u(t)$};
    \node[above,text=teal!75!black] at (10.65,3.75)
      {$v(t-s)$, $t\ge s$};
    \node[above,text=blue!75!black,font=\footnotesize] at (2.6,3.55)
      {follow $u$ up to $s$};
    \node[above,text=teal!75!black,font=\footnotesize] at (8.7,3.45)
      {continue along $v$};
  \end{tikzpicture}
  \caption{Concatenation at a common state.
  The second path is drawn as $v(t-s)$ for $t\ge s$, so its initial
  state is $v(0)=u(s)$. The thick segments select a continuous joined
  path; its admissibility requires a separate concatenation argument.}
  \label{fig:trajectory-concatenation}
\end{figure}
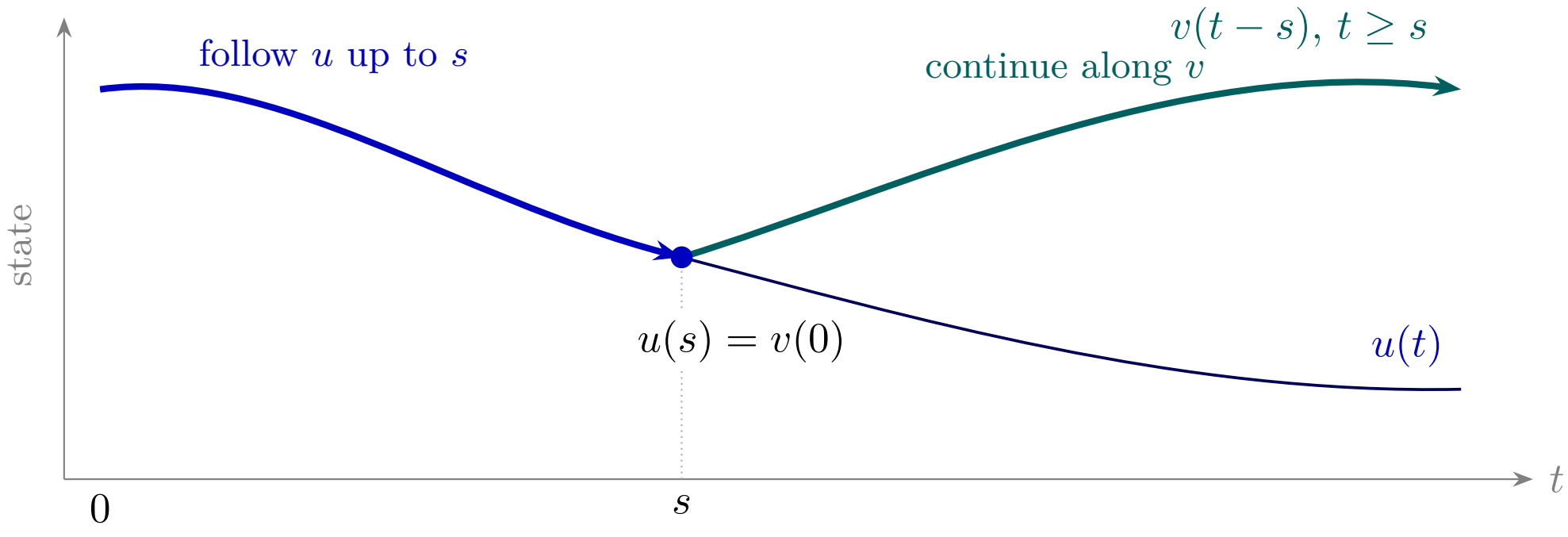

Suppose that two primitive trajectories $u,v$ satisfy $u(s)=v(0)$.
Their smooth approximants $u_j,v_j$ satisfy
\[
 \norm{u_j(s)-v_j(0)}_{H^{-1}}\longrightarrow0,
\]
but the approximating states need not agree for any finite $j$.
Concatenating these approximants may therefore produce a jump at time $s$,
even though the joined limiting path is continuous.
The preceding compactness argument does not guarantee that this path
can be approximated by smooth solutions evolving without restarts.
We therefore allow piecewise restarts and require the total jump size
to tend to zero.
The resulting limits will be called defect trajectories.
The concatenation at a common state is illustrated in
\Cref{fig:trajectory-concatenation}.

\FloatBarrier

\begin{introdefinition} 
\label{def:piecewise-classical}
On $[0,T]$, a piecewise classical approximant consists of a finite
partition $0=t_0<\cdots<t_N=T$, integer cutoffs $\nu_\ell\ge0$, and smooth
solutions $w_\ell$ on the closed intervals $[t_{\ell-1},t_\ell]$ satisfying
\[
 \partial_t w_\ell+\partial_x^3w_\ell
 -3\partial_x(w_\ell^2)+\gamma w_\ell=f_{\nu_\ell},
 \qquad 1\le\ell\le N.
\]
Set $w(t)=w_\ell(t)$ for $t\in[t_{\ell-1},t_\ell)$ and $w(T)=w_N(T)$.
The jump at an interior partition point is
\[
 J_\ell(w)=w_{\ell+1}(t_\ell)-w_\ell(t_\ell),
 \qquad 1\le\ell\le N-1.
\]
Define the total restart error and the minimum force cutoff by
\begin{equation}\label{eq:defect-and-resolution}
 \mathfrak d_T(w)=\sum_{\ell=1}^{N-1}\norm{J_\ell(w)}_{H^{-1}},
 \qquad
 \nu_*(w)=\min_{1\le\ell\le N}\nu_\ell.
\end{equation}
For a continuous path $u:[0,T]\to H^{-1}_0$, set
\[
 \norm{w-u}_{C_*([0,T];H^{-1})}
 =\max_{1\le\ell\le N}
   \sup_{t\in[t_{\ell-1},t_\ell]}\norm{w_\ell(t)-u(t)}_{H^{-1}}.
\]
This comparison includes both one-sided traces at every interior
partition point.
Restrictions to a subinterval retain only pieces whose intersection
with its interior has positive length, using the right trace at the
initial endpoint and the left trace at the terminal endpoint.
\end{introdefinition}

\begin{introdefinition} \label{def:defect-trajectory}
A path $u\in C_{\mathrm{loc}}([0,\infty);H^{-1}_0)$ is a
\emph{defect trajectory} if there exists a sequence of piecewise
classical approximants $w_j$ on $[0,j]$ with the following property.
For every fixed $T>0$, the restrictions $w_j^T=w_j|_{[0,T]}$, defined for
$j\ge T$, satisfy
\begin{equation}\label{eq:defect-convergence}
  \lim_{j\to\infty}\norm{w_j^T-u}_{C_*([0,T];H^{-1})}=0,\qquad
  \lim_{j\to\infty}\mathfrak d_T(w_j^T)=0,\qquad
  \lim_{j\to\infty}\nu_*(w_j^T)=\infty.
\end{equation}
Each approximant has finitely many restarts, but their number need not
be bounded as $j\to\infty$.
The class of all defect trajectories is denoted by $\cK_*^+$.
\end{introdefinition}

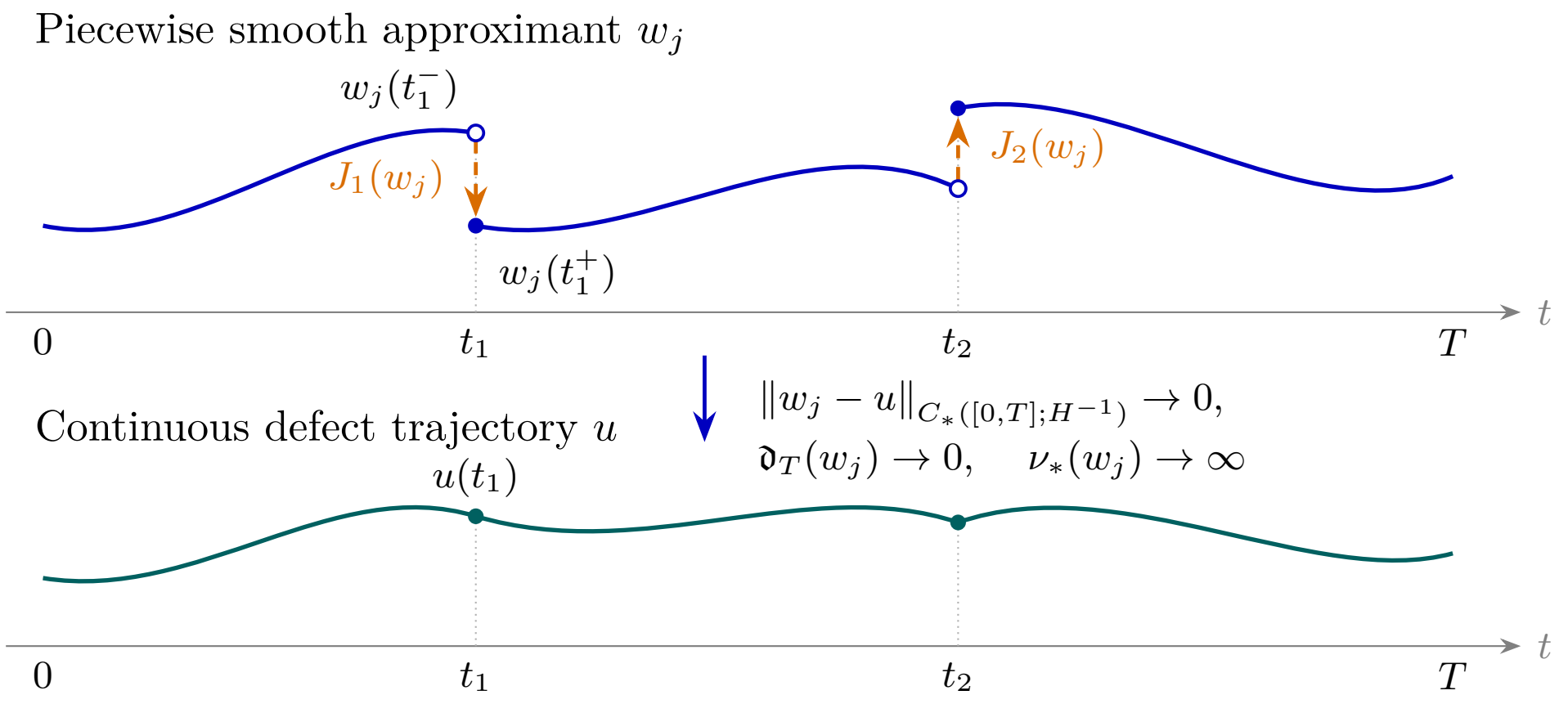
\begin{figure}[htbp]
  \centering
  \begin{tikzpicture}[x=1cm,y=1cm,>=Stealth,
      every node/.style={font=\footnotesize},
      smooth piece/.style={draw=blue!75!black,line width=1pt},
      jump/.style={draw=orange!85!black,densely dashed,line width=.9pt,->}]
    \node[anchor=west,font=\small] at (.45,5.3)
      {Piecewise smooth approximant $w_j$};
    \draw[->,gray] (.35,3.05)--(12.6,3.05) node[right] {$t$};
    \draw[smooth piece] (.65,3.75)
      .. controls (1.85,3.5) and (3.0,4.7) .. (4.15,4.5);
    \draw[smooth piece] (4.15,3.75)
      .. controls (5.45,3.5) and (6.8,4.65) .. (8.05,4.05);
    \draw[smooth piece] (8.05,4.7)
      .. controls (9.3,4.95) and (10.9,3.65) .. (12.05,4.15);
    \draw[gray!60,densely dotted] (4.15,3.05)--(4.15,4.5);
    \draw[gray!60,densely dotted] (8.05,3.05)--(8.05,4.7);
    \draw[jump] (4.15,4.44)--(4.15,3.82)
      node[midway,left=3pt,text=orange!85!black] {$J_1(w_j)$};
    \draw[jump] (8.05,4.12)--(8.05,4.63)
      node[midway,right=3pt,text=orange!85!black] {$J_2(w_j)$};
    \draw[draw=blue!75!black,fill=white,line width=.7pt]
      (4.15,4.5) circle[radius=1.8pt]
      (8.05,4.05) circle[radius=1.8pt];
    \fill[blue!75!black] (4.15,3.75) circle[radius=1.8pt]
      (8.05,4.7) circle[radius=1.8pt];
    \node[above left] at (4.15,4.56) {$w_j(t_1^-)$};
    \node[below right] at (4.2,3.69) {$w_j(t_1^+)$};
    \node[below] at (.65,3.05) {$0$};
    \node[below] at (4.15,3.05) {$t_1$};
    \node[below] at (8.05,3.05) {$t_2$};
    \node[below] at (12.05,3.05) {$T$};

    \draw[->,blue!75!black,line width=.9pt] (6.0,2.7)--(6.0,2.0);
    \node[anchor=west,align=left] at (6.3,2.1)
      {$\norm{w_j-u}_{C_*([0,T];H^{-1})}\to0$,\\[1mm]
       $\mathfrak d_T(w_j)\to0$, \quad $\nu_*(w_j)\to\infty$};
    \node[anchor=west,font=\small] at (.45,2.1)
      {Continuous defect trajectory $u$};
    \draw[->,gray] (.35,.35)--(12.6,.35) node[right] {$t$};
    \draw[teal!75!black,line width=1pt] (.65,.9)
      .. controls (1.9,.7) and (3.0,1.75) .. (4.15,1.4)
      .. controls (5.5,1.0) and (6.85,1.75) .. (8.05,1.35)
      .. controls (9.35,1.8) and (10.95,.8) .. (12.05,1.1);
    \draw[gray!60,densely dotted] (4.15,.35)--(4.15,1.4);
    \draw[gray!60,densely dotted] (8.05,.35)--(8.05,1.35);
    \fill[teal!75!black] (4.15,1.4) circle[radius=1.8pt]
      (8.05,1.35) circle[radius=1.8pt];
    \node[above] at (4.15,1.42) {$u(t_1)$};
    \node[below] at (.65,.35) {$0$};
    \node[below] at (4.15,.35) {$t_1$};
    \node[below] at (8.05,.35) {$t_2$};
    \node[below] at (12.05,.35) {$T$};
  \end{tikzpicture}
  \caption{Piecewise smooth approximation of a continuous defect trajectory.
  Open and filled endpoints denote left and right traces; dashed arrows
  denote jumps. The comparison norm includes both traces at every restart.
  Restart times may depend on $j$.}
  \label{fig:defect-approximation}
\end{figure}

The limiting construction is illustrated in
\Cref{fig:defect-approximation}.
The approximants may jump at restart times, but both traces approach
the same continuous limiting trajectory.
At the same time, the total restart error vanishes and the minimum
force cutoff tends to infinity.

\FloatBarrier

Taking approximants with a single piece gives
$\cK_0^+\subseteq\cK_*^+$.
The defect class is closed under time translation and concatenation.
For the latter property, the additional jump created by joining two
approximants is bounded by their approximation errors at the common
state.
Its size, together with the two original total restart errors, tends
to zero; see Proposition~\ref{prop:defect-structure}.

Moreover, the determinant comparison
\eqref{eq:intro-spectral-comparison} and the dissipation estimate
\eqref{eq:intro-dissipation} remain valid for defect trajectories,
with thresholds $\kappa_R$ uniform for initial states in the
$H^{-1}$ ball of radius $R$.
The additional error across restarts is controlled by the total
$H^{-1}$ jump size and vanishes in the defining limit, see Proposition~\ref{prop:defect-beta}.

We now use defect trajectories to define evolution in the phase space
$H^{-1}_0$:
\begin{equation}\label{eq:multiflow}
 \cS_*(t)u_0=\{u(t):u\in\cK_*^+,\ u(0)=u_0\},
 \qquad t\ge0,\quad u_0\in H^{-1}_0.
\end{equation}
For $B\subset H^{-1}_0$, set
$\cS_*(t)B=\bigcup_{u_0\in B}\cS_*(t)u_0$.
This evolution is called a \emph{strict multivalued semiflow} if all
its values are nonempty and
\begin{equation}\label{eq:strict-semigroup}
 \cS_*(0)u_0=\{u_0\},\qquad
 \cS_*(t+s)u_0=\cS_*(t)\bigl(\cS_*(s)u_0\bigr)
 \quad \text{for} \quad s,t\ge0.
\end{equation}
In the composition law, time translation gives the inclusion from
left to right.
Closure under concatenation gives the reverse inclusion by allowing
any admissible continuation from the intermediate state.

We can now formulate attraction in the phase space. The following definition agrees with
\cite[Definition~2.17]{HanKloeden2023} and
\cite[Theorem~3 and Remark~8]{MelnikValero1998}.
For the associated existence theory for multivalued
attractors, see \cite{MelnikValero1998,MelnikValero2008}.

\begin{introdefinition}
\label{def:global-attractor}
For a strict multivalued semiflow $\cS_*(t)$ on $H^{-1}_0$, a compact
global attractor is a nonempty compact set $\cA_*\subset H^{-1}_0$
that is strictly invariant and attracts every bounded subset of
$H^{-1}_0$.
Thus
\begin{equation}\label{eq:global-attractor-definition}
 \begin{aligned}
 \cS_*(t)\cA_*&=\cA_* &&(t\ge0),\\
 \lim_{t\to\infty}\dist_{H^{-1}}\bigl(\cS_*(t)B,\cA_*\bigr)&=0
 &&\text{for every bounded }B\subset H^{-1}_0,
 \end{aligned}
\end{equation}
where
\[
 \dist_{H^{-1}}(E,F)
 =\sup_{p\in E}\inf_{q\in F}\norm{p-q}_{H^{-1}}
\]
is the one-sided Hausdorff distance. 
\end{introdefinition}

To describe the structure of the attractor, we extend the notion of a
complete trajectory to the defect class.
A path $U\in C(\mathbb R;H^{-1}_0)$ is called a
\emph{complete defect trajectory} if
\[
 U(s+\cdot)|_{[0,\infty)}\in\cK_*^+
 \qquad\text{for every }s\in\mathbb R.
\]
Such a trajectory is \emph{bounded} if
\[
 \sup_{t\in\mathbb R}\norm{U(t)}_{H^{-1}}<\infty.
\]
For paths defined on $\mathbb R$, we also use the notation
$(\mathcal T_sU)(t)=U(t+s)$ for all $s,t\in\mathbb R$.

\begin{theorem}
\label{thm:intro-defect}
Let $\gamma>0$ and $f\in L^2_0(\mathbb T)$.
The evolution \eqref{eq:multiflow} is a strict multivalued semiflow on
$H^{-1}_0$ with closed graph for each fixed $t\ge0$.
It admits a unique compact global attractor $\cA_*$ in the sense of
\Cref{def:global-attractor}.
Moreover,
\[
 \cA_*=\{U(0):U\text{ is a bounded complete defect trajectory}\}.
\]
\end{theorem}

The semiflow and attractor are constructed in
\Cref{thm:strict-multiflow,thm:multivalued-attractor}.
The argument combines the closure properties of the defect class with
the inherited dissipative estimates and time equicontinuity in
Birkhoff coordinates.

The preceding theorem constructs the attractor in $H^{-1}_0$.
We now address its regularity and compare it with the classical
$L^2$ attractor.
For this purpose, write $S_2(t)$ for the classical $L^2_0$ solution
semigroup of \eqref{eq:fdkdv} and $\cA_2$ for its global attractor
(see \Cref{thm:classical-L2}).
A bounded complete trajectory of $S_2$ is a path
$V\in C(\mathbb R;L^2_0)$ such that
\[
 V(t)=S_2(t-s)V(s)\quad \text{for} \quad s\le t ,
 \qquad \text{and}  \qquad
 \sup_{t\in \mathbb R}\norm{V(t)}_{L^2}<\infty.
\]

\begin{theorem} 
\label{thm:intro-L2-identification}
Let $\gamma>0$ and $f\in L^2_0(\mathbb T)$.
Every bounded complete defect trajectory $U$ satisfies
\begin{equation}\label{eq:intro-L2}
 \sup_{t\in\mathbb R}\norm{U(t)}_{L^2}^2
 \le C\gamma^{-2}\norm f_{L^2}^2,
\end{equation}
where $C>0$ is an absolute constant.
For every $u_0\in L^2_0$, the forward defect trajectory from $u_0$ is
unique and coincides with the classical trajectory:
\begin{equation}\label{eq:defect-forward-L2-identification}
 \{u\in\cK_*^+:u(0)=u_0\}=\{S_2(\cdot)u_0\}.
\end{equation}
The bounded complete defect trajectories are precisely the bounded
complete trajectories of $S_2$.
Consequently, $\cA_*=\cA_2$ as subsets of $H^{-1}_0$.
In particular,  $\cA_*=\{0\}$ if $f=0$.
\end{theorem}

The $L^2$ bound follows from \Cref{thm:pullback-L2}.
The identification with the classical dynamics is proved in
Section~\ref{sec:pullback-l2}, using the Birkhoff coordinate equations.

\subsection{Previous results}

The choice of phase space in this paper is motivated by the low-regularity
Cauchy theory for periodic KdV.
Bourgain \cite{Bourgain1993KdV} introduced Fourier restriction spaces
and established global well-posedness in $L^2(\mathbb T)$.
He also treated measure-valued initial data of sufficiently
small norm, combining space-time Fourier estimates with periodic
spectral invariants \cite{Bourgain1997}.
The bilinear estimates of Kenig, Ponce, and Vega
\cite{KenigPonceVega1996} subsequently yielded local well-posedness in
$H^s(\mathbb T)$ for $s\ge-1/2$.
By constructing almost conserved modified energies, Colliander, Keel,
Staffilani, Takaoka, and Tao \cite{CKSTT2003} obtained global
well-posedness in the same range; related periodic multilinear estimates
were developed in \cite{CKSTT2004}.
For the corresponding real-line endpoint theory in
$H^{-3/4}(\R)$, see \cite{Guo2009,Kishimoto2009}.

Further progress at lower regularity relied on the integrable structure.
The spectral approach originates in the inverse-scattering
method of Gardner, Greene, Kruskal, and Miura \cite{Gardner1967}
and the Lax representation \cite{Lax1968}.
For periodic inverse spectral theory and global Birkhoff coordinates,
see \cite{McKeanTrubowitz1976,KappelerPoschel2003}.
Kappeler, M\"ohr, and Topalov \cite{KappelerMoehrTopalov2005}
constructed globally defined real analytic Birkhoff coordinates on
$H^{-1}_0(\mathbb T)$.
Using these coordinates, Kappeler and Topalov
\cite{KappelerTopalov2006} extended the smooth KdV flow to a jointly
continuous group on $H^{-1}(\mathbb T)$.
In Birkhoff coordinates, the conservative dynamics consist of rotations
whose frequencies are determined by the conserved actions.
The analytic properties of the coordinate map and the continuity of the
individual frequencies will be used below to obtain time
equicontinuity on compact subsets of $H^{-1}_0$.

Another approach to the endpoint was developed by Killip and Vi\c{s}an
\cite{KillipVisan2019}, who used commuting Hamiltonian flows to give a
new proof of periodic $H^{-1}$ well-posedness.
Their method builds on the spectral conservation laws introduced by
Killip, Vi\c{s}an, and Zhang \cite{KillipVisanZhang2018} through the
renormalized perturbation determinant.
The determinant bounds in \cite{KillipVisanZhang2018}, together with
the high-frequency compactness arguments in \cite{KillipVisan2019},
provide the spectral tools for our spatial compactness estimates.
Moreover, Molinet \cite{Molinet2012} proved discontinuity of the
periodic KdV flow in $H^s$ for $s<-1$, showing that the endpoint is
sharp in the Sobolev scale.

Turning to the forced and weakly damped periodic KdV equation, Ghidaglia
\cite{Ghidaglia1988,Ghidaglia1994} constructed global attractors and
established finite-dimensional asymptotic behavior at positive Sobolev regularity.
Using Bourgain-space estimates and a decomposition into smoother and
decaying components, Goubet \cite{Goubet2000} proved asymptotic
smoothing in the $L^2$ phase space.
At negative regularity, Tsugawa \cite{Tsugawa2004} used a modified
energy to construct the global attractor in $H^s$, $s>-3/8$, and
identified it with the attractor at nonnegative regularity.
Li and Wu \cite{LiWu2009} combined the $I$-method with an asymptotic
smoothing decomposition to treat $-1/2<s<0$, allowing forcing $f$ in
$H^s_0$.
For $f\in L^2_0$, Yang \cite{Yang2011} reached the endpoint
$s=-1/2$ by combining the $I$-method with a decomposition into a
uniformly $L^2$-bounded part and an exponentially decaying remainder.

A different route to asymptotic compactness is provided by nonlinear
smoothing.
Combining normal-form reductions with Bourgain-space estimates,
Erdo\u{g}an and Tzirakis \cite{ErdoganTzirakis2012} established global
smoothing for periodic KdV.
They also used this approach to construct the $L^2$ global attractor
for the forced and damped KdV equation \cite{ErdoganTzirakis2013}.
In this setting, the solution splits into an exponentially decaying
linear part and a smoother nonlinear remainder, which yields strong
asymptotic compactness.
These results provide the classical $L^2$ attractor used in our
comparison.
Related normal-form methods can be found in
\cite{BabinIlyinTiti2011}.
For spatially analytic forcing, Goubet \cite{Goubet2018} further proved
that the attractor consists of analytic functions.

Related attractor problems have also been studied for other periodic
dispersive equations.
For gKdV equations with suitable polynomial nonlinearities, McConnell
\cite{McConnell2023} combined a gauge transformation, normal-form
analysis, and nonlinear smoothing to construct an $H^1$ global attractor.
For periodic mKdV, Goyal \cite{Goyal2021} used the $I$-method and
trilinear estimates accounting for resonant frequencies to obtain a
global attractor below the energy space.
For the periodic cubic Schr\"odinger equation, Molinet
\cite{Molinet2009Attractor} combined an analysis of weak $L^2$ limits
with Ball's argument to construct an $L^2$ global attractor, and then
used asymptotic smoothing to prove its compactness in $H^2$.
For periodic BBM, Wang and Liu \cite{WangLiu2018BBM} treated
distributional forcing in low-regularity $L^p$-based Sobolev spaces.
Wang \cite{Wang2024BBM} subsequently used a decomposition of the
solution to construct an attractor at sharp low regularity and
quasi-stability estimates to prove its finite fractal dimension.
On $\mathbb T^3$, Bourgain's high--low frequency decomposition also
allowed Wang \cite{Wang2024BBM3D} to construct global attractors below
the energy space.

On the real line, Goubet and Rosa \cite{GoubetRosa2002}
combined the energy equation method with a decomposition of the
solutions to construct a global attractor in $L^2(\R)$ for weakly
damped KdV with $L^2(\R)$ forcing.
The attractor is compact in $H^3(\R)$.
For the weakly damped cubic Schr\"odinger equation with $L^2$ forcing,
Goubet and Molinet
\cite{GoubetMolinet2009} used weak continuity of the $L^2$ flow
and Ball's energy argument to obtain a global attractor in $L^2(\R)$.
The $I$-method and high--low frequency decomposition also allowed
Wang \cite{Wang2014BBM} to construct an $L^2(\R)$ global attractor
for damped BBM with $f\in L^2(\R)$.
For related global attractor results for weakly damped forced gKdV
on the real line, see Wang \cite{Wang2015gKdV}.
 
\subsection{Proof strategy and organization}

Weak damping provides no parabolic regularization,
in contrast to the second-order dissipation in KdV--Burgers
\cite{MolinetVento2011}.
The proof combines spectral dissipation estimates with the Birkhoff coordinate
equations. These tools address different parts of strong path
compactness: the determinant controls the spatial Fourier tails,
whereas the Birkhoff coordinate equations provide time equicontinuity on
compact state sets. We first obtain these estimates for smooth
solutions, then pass to the two trajectory classes.

For the dissipative estimate, we differentiate $\kappa\alpha_\kappa$
along \eqref{eq:fdkdv}. Conservation under the KdV part leaves only
the direction $f-\gamma q$. A first-exit argument ensures that the operator norm of
$A_\kappa(u(t))$ remains sufficiently small.
In this regime, $\alpha_\kappa(u(t))$ is well defined, and
$\kappa\alpha_\kappa(u(t))$ is comparable to
$\norm{u(t)}_{H^{-1}_\kappa}^2$. This gives \Cref{thm:intro-determinant}, uniformly in the force
cutoff. Strong convergence then passes the estimate to primitive
trajectories. For piecewise smooth approximants, the extra term is
bounded by the total $H^{-1}$ restart error, independently of the
number of pieces. It therefore vanishes in the defect limit.

The resulting absorbing bound and Fourier-tail estimate yield
spatial compactness. Indeed, for every bounded initial set
$B\subset H^{-1}_0$,
\[
 \norm{P_{>2\kappa}u(t)}_{H^{-1}}^2
 \le 5\norm{u(t)}_{H^{-1}_\kappa}^2
 \le C\kappa\alpha_\kappa(u(t)).
\]
At fixed sufficiently large $\kappa$, the initial term in
\eqref{eq:intro-dissipation} decays as $t\to\infty$. The remaining
term tends to zero as $\kappa\to\infty$, since
$\norm f_{H^{-1}_\kappa}\to0$. Consequently,
\[
 \lim_{\kappa\to\infty}\limsup_{t\to\infty}
 \sup_{\substack{u\in\cK_*^+\\u(0)\in B}}
 \norm{P_{>2\kappa}u(t)}_{H^{-1}}^2=0.
\]
 Together with the relative compactness of the low-frequency components
in finite-dimensional spaces, this shows that every sequence of late
states has a strongly convergent subsequence in $H^{-1}$.

In order to obtain compactness of paths, we next use the scalar Birkhoff
equations
\[
 \frac{d}{dt}\Omega_n(u(t))
 =i\omega_n(u(t))\Omega_n(u(t))
   +D\Omega_n(u(t))[f-\gamma u(t)].
\]
Here $\Omega_n$ is the $n$th Birkhoff coordinate and $\omega_n$ is its
conservative KdV frequency. On a compact state set, the right-hand
side is uniformly bounded for each fixed $n$. Compactness of the
Birkhoff image controls the remaining coordinates, and continuity of
$\Omega^{-1}$ gives equicontinuity in $H^{-1}$. Arzel\`a--Ascoli theorem then
provides convergence on each fixed time interval. For primitive
trajectories, limits of large time translates give the trajectory
attractor. For the defect class, time translation and concatenation
give the strict multivalued semiflow; the same compactness estimates
yield its global attractor and bounded complete trajectories.
The restart errors also vanish in the scalar coordinate equations and
in the distributional equation, see \Cref{thm:birkhoff-integral,thm:closed-nonlinear-graph} for more details.

Finally, we show that bounded complete trajectories gain $L^2$ regularity by
taking the initial time tend to $-\infty$ at fixed $\kappa$ in the
dissipative estimate. This removes the rough initial-state term and
gives
\[
 \kappa\alpha_\kappa(U(t))
 \le C\gamma^{-2}\norm f_{H^{-1}_\kappa}^2
 \le C\gamma^{-2}\kappa^{-2}\norm f_{L^2}^2.
\]
Combining this estimate with the lower bound in
\eqref{eq:intro-spectral-comparison}, we obtain
\[
 \sum_{|k|\le\kappa}|\widehat U(t,k)|^2
 \le 5\kappa^2\norm{U(t)}_{H^{-1}_\kappa}^2
 \le C\kappa^2\bigl[\kappa\alpha_\kappa(U(t))\bigr]
 \le C\gamma^{-2}\norm{f}_{L^2}^2.
\]
Letting $\kappa\to\infty$ through integers and using the uniformity
in $t$ proves \eqref{eq:intro-L2}.

  To identify these paths with the classical evolution, we combine endpoint analyticity of the Birkhoff map with its one-smoothing property on $L^2$ balls. By utilizing centered differences and interpolation, we get a uniform derivative bound, see Lemma~\ref{lem:Birkhoff-derivative-balls}. After removing the Airy phases, the scalar equations form a locally Lipschitz integral equation in $h^{1/2}$. Its uniqueness identifies the trajectories with the
classical $L^2$ solutions and hence identifies the two attractors.

The paper is organized as follows.
Section~\ref{sec:prelim} records the Birkhoff mapping properties,
perturbation determinant, compactness criteria, and classical $L^2$
dynamics. In section~\ref{sec:determinant}, we establish  dissipation estimates and Fourier-tail estimates. Then, we construct  primitive trajectories and their strong trajectory attractor for \eqref{eq:fdkdv} in Section~\ref{sec:primitive}.  The equations for defect trajectories are derived in  Subsection~\ref{5.1}, thereafter  the strict multivalued semiflow and its compact global attractor are constructed in  Subsection~\ref{5.2}. Finally, Section~\ref{sec:pullback-l2} is devoted to proving $L^2$ regularity and identifying  both the trajectories and the attractor with their classical $L^2$ counterparts.

\section{Preliminaries}\label{sec:prelim}

We first collect the properties of the Birkhoff map and the perturbation
determinant that will be used later. The Birkhoff coordinate equations
will provide time equicontinuity on compact state sets, while the
determinant estimates will yield spatial compactness. We also recall the
compactness criteria and the classical $L^2$ dynamics needed for the
attractor comparison.

\subsection{Birkhoff coordinates}
\label{sec:prelim-coordinates}

The Fourier representation of the Airy flow is
\[
 \widehat u(t,k)=e^{ik^3t}\widehat u(0,k).
\]
Birkhoff coordinates give an analogous rotation formula for the
conservative KdV equation. We write $\Omega(q)=\big(\Omega_1(q), \Omega_2(q), \cdots, \Omega_n(q), \cdots\big)_{n\ge1}$, where
$\Omega_n(q)\in\mathbb C$ denotes the $n$th Birkhoff coordinate of $q$. Along the conservative orbit, the modulus of each coordinate is fixed, and its
angular frequency $\omega_n(q)$ is determined by the conserved actions
$\{|\Omega_n(q)|^2/2\}_{n\ge1}$.

The global Birkhoff map on $H^{-1}_0$ was constructed by Kappeler,
M\"ohr, and Topalov using periodic inverse spectral theory
\cite{KappelerMoehrTopalov2005}. Kappeler and Topalov then extended the
frequencies analytically to $H^{-1}_0$ and obtained a continuous flow
from the rotation formula \cite{KappelerTopalov2006}.
The Riccati representation of distributional potentials
and the Sobolev mapping properties of the Birkhoff map are discussed
in \cite{KappelerTopalovRiccati2005,KappelerSerierTopalov2008}.
For smooth data,
we use the classical Cauchy theory \cite{Kato1979,BonaZhang1996}.
The following statement fixes the normalization used throughout the paper.

\begin{theorem}\label{thm:KT}
 Let $\Phi_t$ denote the KdV flow on $H^{-1}_0(\mathbb T)$.
The family $\{\Phi_t\}_{t\in\mathbb R}$ forms a group, and the map
$(t,q)\mapsto\Phi_t q$ is continuous. There is a real analytic diffeomorphism
\begin{equation}\label{eq:Birkhoff-map}
 \Omega:H^{-1}_0\longrightarrow h^{-1/2}
\end{equation}
such that
\begin{equation}\label{eq:KdV-diagonal}
 \Omega_n(\Phi_tq)=e^{i\omega_n(q)t}\Omega_n(q),
 \qquad q\in H^{-1}_0,\quad t\in\mathbb R,\quad n\ge1.
\end{equation}
For each fixed $n$, the maps $q\mapsto\Omega_n(q)$ and
$q\mapsto\omega_n(q)$ are continuous on $H^{-1}_0$. For smooth $q$,
\begin{equation}\label{eq:Birkhoff-KdV-vector}
 D\Omega_n(q)[-q_{xxx}+6qq_x]=i\omega_n(q)\Omega_n(q).
\end{equation}
\end{theorem}

\begin{proof}
The analytic Birkhoff map is provided by
\cite{KappelerMoehrTopalov2005}, while the continuous KdV group
and the rotation formula are established in
\cite{KappelerTopalov2006}.
After rescaling from the unit circle to $\mathbb T$ and normalizing
the coordinates and frequencies accordingly, these results yield
\eqref{eq:Birkhoff-map} and \eqref{eq:KdV-diagonal}.
 Finally,
differentiating the rotation formula \eqref{eq:KdV-diagonal} at $t=0$ gives
\[
 \begin{aligned}
 D\Omega_n(q)[-q_{xxx}+6qq_x]
  =\left.\frac{\mathrm d}{\mathrm dt}\Omega_n(\Phi_tq)\right|_{t=0} =i\omega_n(q)\Omega_n(q),
 \end{aligned}
\]
which is \eqref{eq:Birkhoff-KdV-vector}.
\end{proof}

The last identity describes the conservative part of the forced
equation in each coordinate. Analyticity also gives, for every compact
$K\subset H^{-1}_0$ and each fixed $n$,
\begin{equation}\label{eq:finite-coordinate-bound}
 \sup_{q\in K}\left(|\Omega_n(q)|+|\omega_n(q)|
       +\norm{D\Omega_n(q)}_{(H^{-1})^*}\right)<\infty.
\end{equation}
These scalar bounds will yield time equicontinuity once spatial
compactness has been established. Later, in
Section~\ref{sec:pullback-l2}, we will assemble the coordinate equations
into an equation in $h^{1/2}$. For that purpose, we need the following
additional properties on $L^2$.

\begin{theorem}\label{thm:L2-Birkhoff-published}
  The Birkhoff map $\Omega$ restricts to a real analytic diffeomorphism
from $L^2_0(\mathbb T)$ onto $h^{1/2}$. Both $\Omega$ and $\Omega^{-1}$ map bounded sets to
bounded sets and are weakly sequentially continuous. Its linearization
at zero is
\[
 \Omega_L=D\Omega(0):L^2_0\longrightarrow h^{1/2},\qquad
 (\Omega_L q)_n=\sqrt{\frac{2}{n}}\,\widehat q(n).
\]
The nonlinear remainder
\begin{equation}\label{eq:L2-one-smoothing}
 \Omega-\Omega_L:L^2_0\longrightarrow h^{3/2}
\end{equation}
is analytic and bounded on bounded sets. Define
\begin{equation}\label{eq:frequency-renormalization}
 \omega_n^*(q)=\omega_n(q)-n^3.
\end{equation}
Then
\begin{equation}\label{eq:L2-frequency-tail}
 \sup_{\norm q_{L^2}\le R}\sup_{n\ge1}n|\omega_n^*(q)|<\infty,
 \qquad R<\infty,
\end{equation}
and $q\mapsto\omega^*(q)$ is locally Lipschitz from $L^2_0$ to
$\ell^\infty$. Moreover, for every norm-compact $K\subset L^2_0$,
\begin{equation}\label{eq:L2-Birkhoff-compact-derivatives}
 \sup_{q\in K}\bigl(\norm{D\Omega(q)}_{\mathcal L(L^2,h^{1/2})}
 +\norm{D^2\Omega(q)}_{\mathcal L^2(L^2,h^{1/2})}\bigr)<\infty.
\end{equation}
\end{theorem}

\begin{proof}
The analytic mapping properties, the formula for $\Omega_L$, and
\eqref{eq:L2-one-smoothing} follow from
\cite[Theorems~1.3--1.4 and Proposition~5.1]{KappelerSchaadTopalov2013},
with the normalization fixed in Theorem~\ref{thm:KT}.
For weak sequential continuity of both $\Omega$ and $\Omega^{-1}$,
see \cite[Corollary~1.3]{KappelerSchaadTopalov2013}.
By \cite[Corollary~1.5(iii)]{KappelerMolnar2018}, we have
\eqref{eq:L2-frequency-tail}. Moreover,
\cite[Theorem~3.6(i)]{KappelerMolnar2018}, with $s=-1/2$ and $p=r=2$,
gives a real analytic map $\omega^*:H^{-1/2}_0\to\ell^2$.
The embeddings $L^2_0\hookrightarrow H^{-1/2}_0$ and
$\ell^2\hookrightarrow\ell^\infty$ then yield the stated local
Lipschitz property.
Finally, \eqref{eq:L2-Birkhoff-compact-derivatives} follows directly
from real analyticity.
\end{proof}

For the forced and damped equation, the Birkhoff coordinate equations
contain the additional term $D\Omega(q)[f-\gamma q]$.
To control this term in $h^{1/2}$ uniformly on bounded subsets of
$L^2_0$, we need a uniform bound for
$D\Omega(q):L^2_0\to h^{1/2}$.
The following lemma derives this bound from the analyticity of
$\Omega:H^{-1}_0\to h^{-1/2}$ and the one-smoothing estimate
\eqref{eq:L2-one-smoothing}.

\begin{lemma}
\label{lem:Birkhoff-derivative-balls}
For every $R>0$, there is $C_R<\infty$ such that
\begin{equation}\label{eq:Birkhoff-derivative-bound}
 \norm{D\Omega(q)[h]}_{h^{1/2}}\le C_R\norm h_{L^2},
 \qquad q,h\in L^2_0,\quad \norm q_{L^2}\le R.
\end{equation}
Consequently,
\[
 \sup_{\norm q_{L^2}\le R}\norm{D\Omega(q)[f-\gamma q]}_{h^{1/2}}
 \le C_R(\norm f_{L^2}+\gamma R).
\]
\end{lemma}

\begin{proof}
Fix $\norm q_{L^2}\le R$ and $\norm h_{L^2}\le1$. Let
\[
 K_R=\{p\in L^2_0:\norm p_{L^2}\le R+1\}.
\]
By weak compactness in $L^2_0$ and the compact embedding
$L^2_0\hookrightarrow H^{-1}_0$, the set $K_R$ is compact
in $H^{-1}_0$.  According to  the endpoint analyticity in   \Cref{thm:KT}, we have
\[
 \sup_{p\in K_R}\norm{D^3\Omega(p)}_{
 \mathcal L^3(H^{-1},h^{-1/2})}<\infty.
\]  

  Since $\Omega_L$ is bounded and linear,
the identity
\[
 D\Omega(q)[h]=\Omega_Lh+D(\Omega-\Omega_L)(q)[h]
\]
reduces the proof to a uniform estimate for $D(\Omega-\Omega_L)(q)[h]$.
For $0<\delta\le1$, let us denote
\[
 d_\delta
 =\frac{\Omega(q+\delta h)-\Omega(q-\delta h)}{2\delta}-\Omega_L(h).
\]
Taylor's formula in $h^{-1/2}$ and the one-smoothing bound \eqref{eq:L2-one-smoothing} give 
\begin{equation}\label{eq:Birkhoff-centered-difference}
 \begin{aligned}
 \norm{d_\delta-D\Omega(q)[h]+\Omega_L(h)}_{h^{-1/2}}
 \le C_R\delta^2,  \qquad
 \norm{d_\delta}_{h^{3/2}} \le C_R\delta^{-1}.
 \end{aligned}
\end{equation}
Here the segment $q+[-\delta,\delta]h$ lies in the $L^2$ ball of
radius $R+1$, and $\norm h_{H^{-1}}\le1$. 

Applying these two estimates
to adjacent dyadic differences and interpolating, we obtain
\[
 \begin{aligned}
 \norm{d_\delta-d_{\delta/2}}_{h^{1/2}}
  \le
 \norm{d_\delta-d_{\delta/2}}_{h^{-1/2}}^{1/2}
 \norm{d_\delta-d_{\delta/2}}_{h^{3/2}}^{1/2} \le C_R\delta^{1/2}.
 \end{aligned}
\]
Since $\norm{d_1}_{h^{1/2}}\le C_R$, it follows that
\[
 \sum_{j\ge0}\norm{d_{2^{-j}}-d_{2^{-j-1}}}_{h^{1/2}}<\infty,
 \qquad
 \sup_{j\ge0}\norm{d_{2^{-j}}}_{h^{1/2}}\le C_R.
\]
Thus $d_{2^{-j}}$ converges in $h^{1/2}$. By the first estimate in
\eqref{eq:Birkhoff-centered-difference}, its limit is
$D\Omega(q)[h]-\Omega_L(h)$, and hence
\[
 \norm{D\Omega(q)[h]-\Omega_L(h)}_{h^{1/2}}\le C_R.
\]
 
By linearity in $h$, this proves
\eqref{eq:Birkhoff-derivative-bound} for every $h\in L^2_0$.
In particular, taking $h=f-\gamma q$ gives
\[
 \norm{D\Omega(q)[f-\gamma q]}_{h^{1/2}}
 \le C_R\norm{f-\gamma q}_{L^2}
 \le C_R\bigl(\norm f_{L^2}+\gamma R\bigr).
\]

This completes the proof.
\end{proof}

\subsection{Perturbation determinant}

We collect the operator estimates and differential formula for the
perturbation determinant introduced in Subsection~\ref{1.2}.
The estimates compare $\kappa\alpha_\kappa(q)$ with
$\norm q_{H^{-1}_\kappa}^2$, providing the spatial control needed
for compactness.  The derivative formula \eqref{eq:alpha-derivative} will be used to
compute the force and damping contributions to
$\frac{d}{dt}\alpha_\kappa(u(t))$ and derive the dissipative estimate;
see the proof of Theorem~\ref{thm:forced-determinant}.
The local derivative bound \eqref{eq:alpha-local-Lipschitz}
controls the change in $\alpha_\kappa$ across each restart
in terms of the $H^{-1}$ norm of the jump; see the proof of
Lemma~\ref{lem:defect-alpha}.

For $B,C\in\HS$, we use the standard Schatten inequalities
\cite{Simon2005}:
\begin{equation}\label{eq:Schatten-22}
 \norm{BC}_{\Ione}\le\norm B_{\HS}\norm C_{\HS},
 \qquad \norm B_{\op}\le\norm B_{\HS}.
\end{equation}
The following periodic Hilbert--Schmidt estimate is given in
\cite[formula~(6.6)]{KillipVisan2019}.

\begin{lemma}
\label{lem:HS-exact}
For every $q\in H^{-1}(\T)$ and $\kappa\ge1$,
\begin{equation}\label{eq:HS-weighted-norm}
 c\kappa^{-1}\norm{q}_{H^{-1}_\kappa}^2
 \le\norm{A_\kappa(q)}_{\HS}^2
 \le C\kappa^{-1}\norm{q}_{H^{-1}_\kappa}^2,
\end{equation}
where $c,C>0$ are absolute constants.
\end{lemma}

  The upper bound in \eqref{eq:HS-weighted-norm} gives
\begin{equation}\label{eq:A-op-small}
 \norm{A_\kappa(q)}_{\op}
 \le \norm{A_\kappa(q)}_{\HS}
 \le C\kappa^{-1/2}\norm q_{H^{-1}_\kappa}
 \le CR\kappa^{-1/2} 
\end{equation}
for $\norm q_{H^{-1}}\le R$. In particular,
\[
 \lim_{\kappa\to\infty}\sup_{\norm q_{H^{-1}}\le R}
 \norm{A_\kappa(q)}_{\op}=0.
\]
Thus, for each $R>0$, there exists $\kappa_R\ge1$, depending only
on $R$, such that the trace series \eqref{trace2aG} converges
for every $\kappa\ge\kappa_R$ and
$\norm q_{H^{-1}}\le R$. Within this region, the determinant is
comparable to its quadratic term. The following lemma gives this
comparison under an operator norm bound; related estimates appear in
\cite[Lemmas~1.4--1.5]{KillipVisanZhang2018}.

\begin{lemma}
\label{lem:det-comparison}
  Let $A\in\HS$ be self-adjoint and satisfy
$\norm A_{\op}\le\eta<1$.   Then
\begin{equation}\label{eq:alpha-HS-comparison}
 c_\eta\norm A_{\HS}^2
 \le \tr(A-\log(I+A))
 \le C_\eta\norm A_{\HS}^2,
\end{equation}
where $c_\eta,C_\eta>0$ depend only on $\eta$.  
\end{lemma}

\begin{proof}
For $|\lambda|\le\eta$, we have
\[
 \lambda-\log(1+\lambda)
 =\lambda^2\int_0^1\frac{s}{1+s\lambda}\dd s.
\]
Hence
\[
 \frac{\lambda^2}{2(1+\eta)}
 \le\lambda-\log(1+\lambda)
 \le\frac{\lambda^2}{2(1-\eta)}.
\]
Applying these bounds to the eigenvalues of $A$ and summing
proves \eqref{eq:alpha-HS-comparison}.
\end{proof}

Combining Lemmas~\ref{lem:HS-exact} and~\ref{lem:det-comparison}, we obtain
\begin{equation}\label{eq:det-weighted-norm}
 c\norm{q}_{H^{-1}_\kappa}^2
 \le \kappa\alpha_\kappa(q)
 \le C\norm{q}_{H^{-1}_\kappa}^2.
\end{equation}
Here $c,C>0$ are absolute once the operator norm is restricted, for
example, to $\norm{A_\kappa(q)}_{\op}\le1/3$. On each bounded
$H^{-1}$ set this restriction holds uniformly for large $\kappa$.

 We now give the derivative formula for $\alpha_\kappa$ and
the local derivative bound mentioned above.

\begin{lemma}\label{lem:alpha-Frechet}
If $\norm{A_\kappa(q)}_{\op}\le\eta<1$, then $\alpha_\kappa$ is
analytic near $q$, and
\begin{equation}\label{eq:alpha-derivative}
 D\alpha_\kappa(q)[h]
 =\tr\bigl(A_\kappa(q)(I+A_\kappa(q))^{-1}A_\kappa(h)\bigr).
\end{equation}
Moreover, whenever $\norm q_{H^{-1}}\le R$ and
$\norm{A_\kappa(q)}_{\op}\le\eta$,
\begin{equation}\label{eq:alpha-local-Lipschitz}
 |D\alpha_\kappa(q)[h]|\le C_{\kappa,\eta,R}\norm h_{H^{-1}}.
\end{equation}
\end{lemma}

\begin{proof}
For every $\ell\ge2$, the Schatten inequalities
\eqref{eq:Schatten-22} give
\[
 \left|\tr\bigl(A_\kappa(q)^{\ell-1}A_\kappa(h)\bigr)\right|
 \le \eta^{\ell-2}
      \norm{A_\kappa(q)}_{\HS}\norm{A_\kappa(h)}_{\HS}.
\]
The trace series can therefore be differentiated locally term by term.
By cyclicity of the trace, we get
\[
 \begin{aligned}
 D\alpha_\kappa(q)[h]
 =\sum_{\ell\ge2}(-1)^\ell
   \tr\bigl(A_\kappa(q)^{\ell-1}A_\kappa(h)\bigr) =\tr\bigl(
   A_\kappa(q)(I+A_\kappa(q))^{-1}A_\kappa(h)
   \bigr),
 \end{aligned}
\]
which proves \eqref{eq:alpha-derivative}.

Finally, \eqref{eq:Schatten-22} and Lemma~\ref{lem:HS-exact} yield
\[
 \begin{aligned}
 |D\alpha_\kappa(q)[h]|
\le \frac{\norm{A_\kappa(q)}_{\HS}
             \norm{A_\kappa(h)}_{\HS}}{1-\eta}
  \le C_{\kappa,\eta,R}\norm h_{H^{-1}},
 \end{aligned}
\]
as required.
\end{proof}

 \begin{remark}For the corresponding differentiation formulas in terms of the
diagonal Green's function, see
\cite[formula~(2.34) and the proof of Proposition~6.4]{KillipVisan2019}.
 \end{remark}

\subsection{Compactness and attractors}

To establish compactness of trajectories, we combine uniform control
of the Fourier tails with equicontinuity in time.
We begin with the periodic Fourier-tail criterion from
\cite[Lemma~4.3 and Section~6]{KillipVisan2019}.

\begin{theorem}
\label{thm:high-energy}
Let $Q\subset H^{-1}_0(\T)$ be bounded. The following conditions are
equivalent:
\begin{align}
 &\lim_{M\to\infty}\sup_{q\in Q}
   \norm{P_{>M}q}_{H^{-1}}=0;
   \label{eq:frequency-tight}\\
 &\lim_{\kappa\to\infty}\sup_{q\in Q}
   \kappa\alpha_\kappa(q)=0.
   \label{eq:beta-tight}
\end{align}
Consequently, $Q$ is relatively compact in $H^{-1}$ if and only if
\eqref{eq:beta-tight} holds.
\end{theorem}

We next turn to equicontinuity in time. For a compact set
$K\subset H^{-1}_0$, by \Cref{thm:KT}, $\Omega(K)$ is compact in $h^{-1/2}$. Given
$\varepsilon>0$, choose a finite $\varepsilon$-net
$z^{(1)},\ldots,z^{(m)}\in\Omega(K)$. For every $z\in\Omega(K)$,
\[
 \begin{aligned}
 \left(\sum_{n>N}\frac{|z_n|^2}{n}\right)^{1/2}
 &\le\min_{1\le\ell\le m}\norm{z-z^{(\ell)}}_{h^{-1/2}}
 +\max_{1\le\ell\le m}
   \left(\sum_{n>N}\frac{|z_n^{(\ell)}|^2}{n}\right)^{1/2}\\
 &\le\varepsilon+\max_{1\le\ell\le m}
   \left(\sum_{n>N}\frac{|z_n^{(\ell)}|^2}{n}\right)^{1/2}.
 \end{aligned}
\]
The maximum tends to zero as $N\to\infty$. Its coordinate tails
therefore satisfy
\[
 a_N(K):=\sup_{q\in K}
 \sum_{n>N}\frac{|\Omega_n(q)|^2}{n},
 \qquad
 \lim_{N\to\infty}a_N(K)=0.
\]
It remains to control finitely many Birkhoff coordinates in time,
as made precise in the following lemma.

\begin{lemma}\label{lem:homeomorphism-equicontinuity}
Let $K\subset H^{-1}_0$ be compact, and let $u_j:[0,T]\to K$.
Assume that, for every $n\ge1$, the family
$\{\Omega_n\circ u_j\}_j$ is equicontinuous on $[0,T]$. Then
\[
 \lim_{\delta\downarrow0}\sup_j
 \sup_{\substack{s,t\in[0,T]\\|t-s|\le\delta}}
 \norm{u_j(t)-u_j(s)}_{H^{-1}}=0.
\]
\end{lemma}

\begin{proof}
For every $N\ge1$ and $s,t\in[0,T]$, it is easy to see that
\[
 \begin{aligned}
 \norm{\Omega(u_j(t))-\Omega(u_j(s))}_{h^{-1/2}}^2
 &\le
 \sum_{n=1}^N
 \frac{|\Omega_n(u_j(t))-\Omega_n(u_j(s))|^2}{n}
 +4a_N(K).
 \end{aligned}
\]
Given $\varepsilon>0$, choose $N$ such that
$4a_N(K)<\varepsilon^2/2$.
Equicontinuity of the first $N$ coordinates then gives
$\delta>0$ such that
\[
 \sup_j
 \sup_{\substack{s,t\in[0,T]\\|t-s|\le\delta}}
 \sum_{n=1}^N
 \frac{|\Omega_n(u_j(t))-\Omega_n(u_j(s))|^2}{n}
 <\frac{\varepsilon^2}{2}.
\]
Consequently,
\[
 \sup_j
 \sup_{\substack{s,t\in[0,T]\\|t-s|\le\delta}}
 \norm{\Omega(u_j(t))-\Omega(u_j(s))}_{h^{-1/2}}
 <\varepsilon.
\]
Thus $\{\Omega\circ u_j\}_j$ is equicontinuous in $h^{-1/2}$.
Uniform continuity of $\Omega^{-1}$ on the compact set $\Omega(K)$
gives the asserted equicontinuity in $H^{-1}$.
\end{proof}

We now turn to global attractors in the phase space and recall
the properties used in their construction. A single-valued semigroup can be viewed as a multivalued evolution whose images consist of a single state.
Thus \Cref{def:global-attractor} applies to both settings.
Attraction is measured in the norm of the phase space:
$H^{-1}_0$ for the generalized dynamics considered here and
$L^2_0$ for the classical dynamics below.

The construction below uses absorbing sets and asymptotic
compactness.
General accounts of these notions and their role in
attractor theory can be found in \cite{Hale1988,Temam1997,Robinson2001};
see also the survey \cite{Zelik2023}.
A bounded set $B_0\subset X$ is called absorbing
if, for every bounded set $B\subset X$, there exists $t_B\ge0$
such that
\[
 S(t)B\subset B_0
 \qquad\text{for all }t\ge t_B.
\]
The evolution $S(t)$ is called asymptotically compact if,
whenever
\[
 \{x_j\}\subset X\text{ is bounded},\qquad
 t_j\to\infty,\qquad y_j\in S(t_j)x_j,
\]
there exist a subsequence $\{y_{j_k}\}$ and a state $y\in X$
such that
\[
 y_{j_k}\longrightarrow y\quad\text{in }X.
\]

We next recall the classical $L^2_0$ dynamics and its global
attractor, which will later be compared with the attractor
constructed in $H^{-1}_0$. The following theorem summarizes the properties needed below; see \cite{Bourgain1993KdV,ErdoganTzirakis2013}.

\begin{theorem}\label{thm:classical-L2}
Let $\gamma>0$ and $f\in L^2_0(\T)$. Equation~\eqref{eq:fdkdv}
defines a continuous semigroup $S_2(t)$ on $L^2_0(\T)$.
If smooth mean-zero initial states $q_j$ and forces $f_j$ satisfy
\[
 \lim_{j\to \infty} \|q_j-q\|_{L^2_0}=0,\qquad  \text{and} \qquad
 \lim_{j\to \infty} \|f_j-f\|_{L^2_0}=0,
\]
then the corresponding smooth solutions satisfy
\[
 u_j\longrightarrow S_2(\cdot)q
 \quad\text{in }C([0,T];L^2_0) 
\]
for  $0<T<\infty$. The semigroup has a compact global attractor $\cA_2$, and
\begin{equation}\label{eq:L2-attractor-kernel}
 \cA_2
 =\{V(0):V\text{ is a bounded complete }S_2\text{-trajectory}\}.
\end{equation}
\end{theorem}

\section{Dissipative estimates and compactness}\label{sec:determinant}\label{sec:bootstrap}

In this section, we derive a dissipative estimate for
$\alpha_\kappa$, from which an absorbing bound and uniform
high-frequency tail estimates will follow. These estimates
provide the spatial control needed for the trajectory
construction in Section~\ref{sec:compactness}.

We begin with the following conservation law for the periodic
KdV flow; see \cite[Corollary~6.5]{KillipVisan2019}.
This property ensures that, along the forced--damped evolution,
only the force and damping contribute to the time derivative
of $\alpha_\kappa$.

\begin{proposition}
\label{lem:kdv-alpha-conservation}
Let $q\in C^1(I;C^\infty(\T))$ solve the KdV equation
\begin{equation}\label{eq:unforced-kdv}
 q_t=-q'''+6qq'.
\end{equation}
If $\norm{A_\kappa(q(t))}_{\op}<1$ on $I$, then
\begin{equation}\label{eq:alpha-conserved}
 \frac{\dd}{\dd t}\alpha_\kappa(q(t))=0,
 \qquad t\in I.
\end{equation}
In particular, \eqref{eq:alpha-conserved} holds along every smooth
 KdV orbit for all sufficiently large $\kappa$.
\end{proposition}

We now estimate the force and damping contributions to
$\frac{\dd}{\dd t}\alpha_\kappa(u(t))$ as long as
$\norm{A_\kappa(u(t))}_{\op}$ is sufficiently small.
For sufficiently large $\kappa$, the resulting bound ensures
that this smallness condition is preserved for all $t\ge0$.
The dissipative estimate therefore holds globally in time,
uniformly in the Fourier cutoff parameter $\nu$ of the forcing.

\begin{theorem}
\label{thm:forced-determinant}
Let $\gamma>0$, let $g\in C^\infty(\T;\mathbb R)\cap L^2_0$, and let $u$
be a smooth real-valued solution with zero spatial mean
on a time interval $I$ of
\begin{equation}\label{eq:smooth-forced}
 u_t=-u_{xxx}+6uu_x-\gamma u+g.
\end{equation}
Then the following statements hold.

\begin{enumerate}[label=\textup{(\roman*)},leftmargin=2.6em]

\item Fix $\kappa\ge1$, and let $J\subset I$ be an interval.
Assume that, 
$
 \norm{A_\kappa(u(t))}_{\op}<1
$ for every $t\in J$. Then, we have
\begin{equation}\label{eq:alpha-diff}
 \frac{\dd}{\dd t}\alpha_\kappa(u(t))
 =\tr\bigl(
 A_\kappa(u(t))(I+A_\kappa(u(t)))^{-1}
 A_\kappa(-\gamma u(t)+g)
 \bigr).
\end{equation}

\item Fix $\kappa\ge1$ and $0<\eta<1$, and let $J\subset I$
be an interval. Assume that
$
 \norm{A_\kappa(u(t))}_{\op}\le\eta
$ for every $t\in J$.
Then, 
\begin{equation}\label{eq:alpha-coercive}
 \frac{\dd}{\dd t}\alpha_\kappa(u(t))
 +c_\eta\gamma\alpha_\kappa(u(t))
 \le C_\eta\gamma^{-1}\norm{A_\kappa(g)}_{\HS}^2,
\end{equation}
where $c_\eta,C_\eta>0$ depend only on $\eta$.

\item Let $R>0$ and fix $f\in L^2_0(\T)$.
Assume that $[0,\infty)\subset I$, $\norm{u(0)}_{H^{-1}}\le R$, 
and  
\begin{equation}\label{eq:force-family-tail}
 \norm{g}_{H^{-1}_\kappa}^2
 \le\norm{f}_{H^{-1}_\kappa}^2
 \qquad\text{for every }\kappa\ge1.
\end{equation}
Then there exists
$
 \kappa_*=\kappa_*(R,\gamma,\norm f_{H^{-1}})\ge1
$
such that, for every $\kappa\ge\kappa_*$,
\[
 \norm{A_\kappa(u(t))}_{\op}\le\frac14
 \qquad\text{for all }t\ge0.
\]
For these values of $\kappa$, $\alpha_\kappa(u(t))$ is
well defined  and satisfies
\begin{equation}\label{eq:global-beta}
\begin{aligned}
 \kappa\alpha_\kappa(u(t))
  \le e^{-c_0\gamma(t-s)}
       \kappa\alpha_\kappa(u(s)) +C_0\gamma^{-2}\norm{f}_{H^{-1}_\kappa}^2
       \bigl(1-e^{-c_0\gamma(t-s)}\bigr)
\end{aligned}
\end{equation}
for $0\le s\le t$.
Here $c_0,C_0>0$ are absolute constants.
The threshold $\kappa_*$ and these constants are uniform over
all $g$ satisfying \eqref{eq:force-family-tail}.
In particular, they are independent of $\nu$ when
$g=P_{\le\nu}f$.

\end{enumerate}
\end{theorem}

\begin{proof}
\textup{(i)}
It follows from Proposition~\ref{lem:kdv-alpha-conservation}
that
\[
 D\alpha_\kappa(u)[-u_{xxx}+6uu_x]=0.
\]
Thus, by the chain rule and \eqref{eq:alpha-derivative}, we have
\[
 \begin{aligned}
 \frac{\dd}{\dd t}\alpha_\kappa(u)
  =D\alpha_\kappa(u)[-\gamma u+g] =\tr\bigl(
 A_\kappa(u)(I+A_\kappa(u))^{-1}
 A_\kappa(-\gamma u+g)
 \bigr) 
 \end{aligned}
\]
which proves \eqref{eq:alpha-diff}.

\textup{(ii)}
From the spectral theorem, one has
\begin{equation}\label{eq:damping-trace-lower}
 \tr\bigl(A_\kappa(u)^2(I+A_\kappa(u))^{-1}\bigr)
 \ge\frac{1}{1+\eta}\norm{A_\kappa(u)}_{\HS}^2.
\end{equation}
For the force term,   the Schatten inequalities and Young's inequality yield
\begin{equation}\label{eq:force-trace-upper}
 \begin{aligned}
  \abs{\tr\bigl(
 A_\kappa(u(t))(I+A_\kappa(u(t)))^{-1}A_\kappa(g)
 \bigr)} 
 & \le
 \frac{\norm{A_\kappa(u(t))}_{\HS}
       \norm{A_\kappa(g)}_{\HS}}{1-\eta}\\
 & \le
 \frac{\gamma}{2(1+\eta)}
 \norm{A_\kappa(u)}_{\HS}^2
 +C_\eta\gamma^{-1}\norm{A_\kappa(g)}_{\HS}^2.
 \end{aligned}
\end{equation}
 
Substituting these estimates into \eqref{eq:alpha-diff}
and applying the upper bound in
\eqref{eq:alpha-HS-comparison}, we obtain
\eqref{eq:alpha-coercive}.

\textup{(iii)}
By \eqref{eq:HS-weighted-norm}, the initial operator norm
is less than $1/4$ for sufficiently large $\kappa$.
On any interval starting at zero where
$\norm{A_\kappa(u)}_{\op}\le1/4$, we have
\[
 \begin{aligned}
 \norm{A_\kappa(u)}_{\op}^2  \leq C 
 \alpha_\kappa(u(t))
 &\le C\left(
 \alpha_\kappa(u(0))
 +\gamma^{-2}\norm{A_\kappa(g)}_{\HS}^2
 \right)\\
 &\le\frac{C}{\kappa}
 \left(R^2+\gamma^{-2}\norm f_{H^{-1}}^2\right).
 \end{aligned}
\]
The second inequality follows from part \textup{(ii)}
with $\eta=1/4$, Gronwall's inequality, and
\eqref{eq:alpha-HS-comparison}.
For the last, we used \eqref{eq:HS-weighted-norm}
and \eqref{eq:force-family-tail}.

Now choose
$\kappa_*=\kappa_*(R,\gamma,\norm f_{H^{-1}})$
sufficiently large that the last expression is at most
$1/64$ whenever $\kappa\ge\kappa_*$.
This improves the assumed operator norm bound from $1/4$
to $1/8$. A continuity argument then shows that
$\norm{A_\kappa(u)}_{\op}<1/4$ for all $t\ge0$. 

Applying Gronwall's inequality to
\eqref{eq:alpha-coercive} with $\eta=1/4$, we obtain
\[
 \begin{aligned}
 \alpha_\kappa(u(t))
 \le e^{-c_0\gamma(t-s)}\alpha_\kappa(u(s)) +C\gamma^{-2}\norm{A_\kappa(g)}_{\HS}^2
 \bigl(1-e^{-c_0\gamma(t-s)}\bigr).
 \end{aligned}
\]
Multiplying this estimate by $\kappa$ and using
\[
 \kappa\norm{A_\kappa(g)}_{\HS}^2
 \le C\norm g_{H^{-1}_\kappa}^2
 \le C\norm f_{H^{-1}_\kappa}^2
\]
gives \eqref{eq:global-beta}.
All constants are absolute because $\eta=1/4$ is fixed.
Their uniformity, as well as that of $\kappa_*$, follows
from \eqref{eq:force-family-tail}.
In particular, the estimates are uniform in $\nu$
when $g=P_{\le\nu}f$.
\end{proof}

\begin{corollary}\label{cor:absorbing}
For every $R>0$ there are a positive constant $C_R<\infty$ and absolute $c,C>0$ such that all
smooth solutions covered by \Cref{thm:forced-determinant}\textup{(iii)}
satisfy
\begin{equation}\label{eq:Hminus1-absorbing}
 \norm{u(t)}_{H^{-1}}^2
 \le C_R e^{-c\gamma t}
   +C\gamma^{-2}\norm f_{L^2}^2,
 \qquad t\ge0.
\end{equation}
The constant $C_R$ may depend on $R$, $\gamma$, and
$\norm f_{H^{-1}}$.
Consequently, there is a bounded ball $B_*\subset H^{-1}_0$ that absorbs
every bounded initial set. Its radius depends only on $\gamma$ and
$\norm f_{L^2}$.
\end{corollary}

\begin{proof}
Fix $\kappa_0\ge\kappa_*(R,\gamma,\norm f_{H^{-1}})$, 
 it is easy to see that
\begin{equation}\label{eq:Hminus1-from-fixed-weight}
 \norm q_{H^{-1}}^2
 \le4\kappa_0^2\norm q_{H^{-1}_{\kappa_0}}^2 
\end{equation}
and
\[
 \kappa_0^2\norm f_{H^{-1}_{\kappa_0}}^2
 =\sum_k\frac{\kappa_0^2}{k^2+4\kappa_0^2}
       |\widehat f(k)|^2
 \le\frac14\norm f_{L^2}^2.
\]
Using these bounds together with
\eqref{eq:det-weighted-norm} and \eqref{eq:global-beta},
we obtain
\[
 \begin{aligned}
 \norm{u(t)}_{H^{-1}}^2
    \le C\kappa_0^2
       \norm{u(t)}_{H^{-1}_{\kappa_0}}^2 
& \le C\kappa_0^2e^{-c\gamma t}
       \norm{u(0)}_{H^{-1}_{\kappa_0}}^2
    +C\gamma^{-2}\kappa_0^2
       \norm f_{H^{-1}_{\kappa_0}}^2\\
 &\le C\kappa_0^2R^2e^{-c\gamma t}
      +C\gamma^{-2}\norm f_{L^2}^2.
 \end{aligned}
\]
which proves \eqref{eq:Hminus1-absorbing}.
The first term tends to zero as $t\to\infty$, while the
last term is independent of $R$.
Hence every ball centered at zero with squared radius
strictly larger than $C\gamma^{-2}\norm f_{L^2}^2$
is absorbing.
\end{proof}

The following corollary provides the high-frequency tail
control needed to establish spatial compactness.

\begin{corollary}
\label{cor:asymptotic-tail}
Let $\gamma>0$, $f\in L^2_0(\T)$, and $R>0$.
In the suprema below, $u$ ranges over all global smooth
real-valued solutions of \eqref{eq:smooth-forced} with
$\norm{u(0)}_{H^{-1}}\le R$ and smooth mean-zero forces
$g$ satisfying \eqref{eq:force-family-tail}.
Then, for every
$\kappa\ge\kappa_*(R,\gamma,\norm f_{H^{-1}})$, we have
\begin{equation}\label{eq:asymptotic-beta-source}
 \limsup_{t\to\infty}\sup_u
 \kappa\alpha_\kappa(u(t))
 \le C\gamma^{-2}\norm f_{H^{-1}_\kappa}^2,
\end{equation}
where $C>0$ is an absolute constant.
In particular,
\begin{equation}\label{eq:asymptotic-beta-tail}
 \lim_{\kappa\to\infty}\limsup_{t\to\infty}
 \sup_u\kappa\alpha_\kappa(u(t))=0.
\end{equation}
\end{corollary}

\begin{proof}
Taking $s=0$ in \eqref{eq:global-beta} and using
\eqref{eq:det-weighted-norm}, we get
\[
 \sup_u\kappa\alpha_\kappa(u(t))
 \le CR^2e^{-c_0\gamma t}
     +C_0\gamma^{-2}\norm f_{H^{-1}_\kappa}^2 
\]
for every
$\kappa\ge\kappa_*(R,\gamma,\norm f_{H^{-1}})$.
Letting $t\to\infty$ with $\kappa$ fixed proves
\eqref{eq:asymptotic-beta-source}.
Since
\[
 \norm f_{H^{-1}_\kappa}^2
 \le\frac{1}{4\kappa^2}\norm f_{L^2}^2,
\]
we may now let $\kappa\to\infty$ to obtain
\eqref{eq:asymptotic-beta-tail}.
\end{proof}

\section{Primitive trajectories and the trajectory attractor}
\label{sec:compactness}\label{sec:primitive}

In this section, we construct primitive trajectories for every
initial state in $H^{-1}_0$ and show that the class $\cK_0^+$
admits a strong trajectory attractor.

\subsection{Primitive trajectories}

To construct the limits required in Definition~\ref{def:primitive},
we first establish compactness of smooth approximating solutions
on finite time intervals.

Fix $f\in L^2_0(\T)$ and $\gamma>0$.
For each integer $\nu\ge1$ and smooth mean-zero initial datum
$q_\nu$, let $u_\nu$ solve
\begin{equation}\label{eq:nu-problem}
 \begin{cases}
 (u_\nu)_t+(u_\nu)_{xxx}-3\partial_x(u_\nu^2)
 +\gamma u_\nu=f_\nu, & t>0,\\
 u_\nu(0)=q_\nu,
 \end{cases}
\end{equation}
where $f_\nu=P_{\le\nu}f$.
For each fixed $\nu$, classical KdV energy estimates ensure
that $u_\nu$ is smooth and defined for all $t\ge0$.

We first consider compactness in space.
When the initial states range over a compact subset of
$H^{-1}_0$, their high-frequency tails are uniformly small.
The following proposition shows that this property is
preserved on finite time intervals, uniformly in $\nu$.

\begin{proposition}
\label{prop:finite-spatial-tightness}
Let $Q\subset H^{-1}_0(\T)$ be compact, and let
$q_\nu\in Q\cap C^\infty(\T)$.
Let $u_\nu$ be the solution of \eqref{eq:nu-problem}.
Then, for every $T>0$, we have
\begin{equation}\label{eq:finite-window-tail}
 \lim_{\kappa\to\infty}\sup_\nu\sup_{0\le t\le T}
 \kappa\alpha_\kappa(u_\nu(t))=0.
\end{equation}
Moreover, the set
\begin{equation}\label{eq:state-hull}
 K_T^0=\{u_\nu(t):\nu\in\N,\ 0\le t\le T\}
 \subset H^{-1}_0
\end{equation}
is relatively compact in $H^{-1}_0$.  
\end{proposition}

\begin{proof}
Since $Q$ is compact, \Cref{thm:high-energy} implies
\begin{equation}\label{eq:Q-beta-tail}
 \lim_{\kappa\to\infty}
 \sup_{q\in Q}\kappa\alpha_\kappa(q)=0.
\end{equation}
For sufficiently large $\kappa$, \eqref{eq:global-beta}
applies uniformly in $\nu$ and gives
\begin{equation}\label{eq:finite-beta-bound}
 \kappa\alpha_\kappa(u_\nu(t))
 \le e^{-c_0\gamma t}\kappa\alpha_\kappa(q_\nu)
     +C\gamma^{-2}\norm f_{H^{-1}_\kappa}^2,
 \qquad 0\le t\le T.
\end{equation}
Taking suprema and using $q_\nu\in Q$, we get
\[
 \begin{aligned}
 \sup_\nu\sup_{0\le t\le T}
 \kappa\alpha_\kappa(u_\nu(t))
 &\le\sup_{q\in Q}\kappa\alpha_\kappa(q)
     +C\gamma^{-2}\norm f_{H^{-1}_\kappa}^2\\
 &\le\sup_{q\in Q}\kappa\alpha_\kappa(q)
     +\frac{C}{\gamma^2\kappa^2}\norm f_{L^2}^2.
 \end{aligned}
\]
Letting $\kappa\to\infty$ and using
\eqref{eq:Q-beta-tail} proves \eqref{eq:finite-window-tail}.

The absorbing estimate \eqref{eq:Hminus1-absorbing}
also gives
\[
 \sup_\nu\sup_{0\le t\le T}
 \norm{u_\nu(t)}_{H^{-1}}<\infty.
\]
Together with \eqref{eq:finite-window-tail}, this allows
us to apply \Cref{thm:high-energy} and conclude that
$K_T^0$ is relatively compact in $H^{-1}_0$.
\end{proof}

Having established spatial compactness, we now turn to
equicontinuity in time. By
Lemma~\ref{lem:homeomorphism-equicontinuity}, it suffices
to prove equicontinuity of each fixed Birkhoff coordinate,
which we obtain from its evolution equation.

\begin{proposition}
\label{prop:Birkhoff-equicontinuity}
Under the hypotheses of Proposition~\ref{prop:finite-spatial-tightness},
\[
 \lim_{\delta\downarrow0}\sup_\nu
 \sup_{\substack{s,t\in[0,T]\\|t-s|\le\delta}}
 \norm{u_\nu(t)-u_\nu(s)}_{H^{-1}}=0.
\]
\end{proposition}

\begin{proof}
Let $K_T=\overline{K_T^0}$, where the closure is taken
in $H^{-1}$.
By Proposition~\ref{prop:finite-spatial-tightness},
$K_T$ is compact.
For each fixed $n\ge1$, the chain rule and
\eqref{eq:Birkhoff-KdV-vector} give
\begin{equation}\label{eq:smooth-Birkhoff-ODE}
 \frac{\dd}{\dd t}\Omega_n(u_\nu)
 =i\omega_n(u_\nu)\Omega_n(u_\nu)
  +D\Omega_n(u_\nu)[f_\nu-\gamma u_\nu].
\end{equation}

Since $K_T$ is compact, $\omega_n\Omega_n$ and
$D\Omega_n$ are bounded on $K_T$.
Note that
\[
 \sup_\nu\norm{f_\nu}_{H^{-1}}
 \le\norm f_{H^{-1}}.
\]
Hence, for each fixed $n\ge1$, there exists $C_{n,T}>0$
such that
\[
 \left|\frac{\dd}{\dd t}\Omega_n(u_\nu)\right|
 \le C_{n,T},
 \qquad \nu\in\N,\quad t\in[0,T].
\]
Integrating, we obtain
\begin{equation}\label{eq:coordinate-Lipschitz}
 |\Omega_n(u_\nu(t))-\Omega_n(u_\nu(s))|
 \le C_{n,T}|t-s|,
 \qquad s,t\in[0,T],
\end{equation}
where $C_{n,T}$ is independent of $\nu$.

Thus each fixed coordinate family is equicontinuous.
As all $u_\nu$ take values in $K_T$,
Lemma~\ref{lem:homeomorphism-equicontinuity}
gives the required equicontinuity in $H^{-1}$.
\end{proof}

The preceding two propositions give spatial compactness and
equicontinuity in time, so we can now apply the
Arzel\`a--Ascoli theorem.

\begin{theorem}\label{prop:finite-compactness}
Let $Q\subset H^{-1}_0(\T)$ be compact, and let
$q_\nu\in Q\cap C^\infty(\T)$ for every $\nu\in\N$.
Let $u_\nu$ be the corresponding solution of
\eqref{eq:nu-problem}.
Then, for every $T>0$, the family
\begin{equation}\label{eq:finite-time-compactness}
 \{u_\nu|_{[0,T]}:\nu\in\N\}
\end{equation}
is relatively compact in $C([0,T];H^{-1}_0)$.
\end{theorem}

\begin{proof}
Fix $T>0$.
By Proposition~\ref{prop:finite-spatial-tightness}, the set
\[
 \{u_\nu(t):\nu\in\N,\ 0\le t\le T\}
\]
is relatively compact in $H^{-1}_0$.
Moreover, Proposition~\ref{prop:Birkhoff-equicontinuity}
ensures that the family $\{u_\nu\}$ is equicontinuous
on $[0,T]$.
The conclusion now follows from the Arzel\`a--Ascoli theorem.
\end{proof}

We now use this compactness result to construct primitive
trajectories for arbitrary initial states in $H^{-1}_0$.

\begin{proposition}
\label{prop:primitive-existence}
For every $u_0\in H^{-1}_0(\T)$, there exists a primitive
trajectory $u\in\cK_0^+$ such that $u(0)=u_0$.
\end{proposition}

\begin{proof}
Set $q_\nu=P_{\le\nu}u_0$, and let $u_\nu$ solve
\eqref{eq:nu-problem}.
Since $q_\nu\to u_0$ in $H^{-1}$, the set
$\{q_\nu:\nu\ge1\}\cup\{u_0\}$ is compact in $H^{-1}_0$.

Applying Theorem~\ref{prop:finite-compactness} on
$[0,1],[0,2],\ldots$ and taking a diagonal subsequence,
we obtain a path
$u\in C_{\mathrm{loc}}([0,\infty);H^{-1}_0)$ such that
\[
 \lim_{\ell\to\infty}
 \sup_{0\le t\le T}
 \norm{u_{\nu_\ell}(t)-u(t)}_{H^{-1}}=0
 \qquad\text{for every }T>0.
\]
In particular,
$u(0)=\lim_{\ell\to\infty}q_{\nu_\ell}=u_0$.
Definition~\ref{def:primitive} now gives $u\in\cK_0^+$.
\end{proof}

We next establish the dissipative estimates and Birkhoff
coordinate identities needed to study the long-time behavior
of primitive trajectories.

\begin{proposition}
\label{prop:primitive-inherited}
Let $R>0$, and let $u\in\cK_0^+$ satisfy
$\norm{u(0)}_{H^{-1}}\le R$.
For all sufficiently large $\kappa$, the determinant is
defined along $u$, and
\eqref{eq:intro-spectral-comparison} and
\eqref{eq:global-beta} hold.
The absorbing bound \eqref{eq:Hminus1-absorbing} also holds.
The threshold $\kappa_{*}$ and constants can be chosen uniformly over
all such trajectories.

Moreover, for every $n\ge1$ and $0\le s\le t$, we have
\begin{equation}\label{eq:primitive-Birkhoff-integral}
 \Omega_n(u(t))-\Omega_n(u(s))
 =\int_s^t
 \left[
 i\omega_n(u(\tau))\Omega_n(u(\tau))
 +D\Omega_n(u(\tau))[f-\gamma u(\tau)]
 \right]\dd\tau.
\end{equation}
\end{proposition}

\begin{proof}
Fix $T>0$, and choose an  approximating sequence $\{u_j\}$ for $u$.
Then
\[
 \lim_{j\to\infty}
 \norm{u_j-u}_{C([0,T];H^{-1})}=0,
 \qquad
 \lim_{j\to\infty}
 \norm{f_{\nu_j}-f}_{H^{-1}}=0.
\]
After discarding finitely many terms, we may assume that
$\norm{u_j(0)}_{H^{-1}}\le R+1$.
Uniform convergence of the paths shows that
$u([0,T])\cup\bigcup_{j\ge1}u_j([0,T])$
has a compact closure in $H^{-1}$, which we denote by $K_T$.

For $\kappa\ge\kappa_*(R+1,\gamma,\norm f_{H^{-1}})$,
Theorem~\ref{thm:forced-determinant}\textup{(iii)}
and \eqref{eq:HS-weighted-norm} imply
\[
 \norm{A_\kappa(q)}_{\op}\le\frac14
 \qquad\text{for every }q\in K_T.
\]
Thus $\alpha_\kappa$ is defined and uniformly continuous
on $K_T$, so
\[
 \lim_{j\to\infty}
 \norm{\alpha_\kappa(u_j)-\alpha_\kappa(u)}_{C([0,T])}=0.
\]
The weighted Sobolev norms are also continuous in $H^{-1}$.
We can therefore pass to the limit in the smooth estimates
to obtain \eqref{eq:intro-spectral-comparison} and
\eqref{eq:global-beta}.
The absorbing bound \eqref{eq:Hminus1-absorbing}
passes to the limit in the same way.
All these estimates are uniform for
$\norm{u(0)}_{H^{-1}}\le R$, since the approximating
initial states satisfy the common bound $R+1$.

Now fix $n\ge1$.
The maps $\Omega_n$, $\omega_n$, and $D\Omega_n$ are
uniformly continuous on $K_T$, with operator norm used
for $D\Omega_n$.
The convergence of the paths and forces therefore gives
uniform convergence of $\Omega_n(u_j)$ and of the
right-hand sides of \eqref{eq:smooth-Birkhoff-ODE}
on $[0,T]$.
Integrating this equation over $[s,t]$ and passing to
the limit, we obtain \eqref{eq:primitive-Birkhoff-integral}
for $0\le s\le t\le T$.
Since $T>0$ was arbitrary, the estimates and the identity
hold on $[0,\infty)$.
\end{proof}

Finally, we establish the closedness and time-translation
properties needed for the construction of the trajectory
attractor.

\begin{proposition}
\label{prop:primitive-closed}
The class $\cK_0^+$ is sequentially closed in
$C_{\mathrm{loc}}([0,\infty);H^{-1}_0)$.
Moreover, it is preserved under forward time translations:
\[
 \mathcal T_s\cK_0^+\subset\cK_0^+,
 \qquad s\ge0.
\]
\end{proposition}

\begin{proof}
We first prove sequential closedness.
Let $\{u^m\}\subset\cK_0^+$ converge to $u$ in
$C_{\mathrm{loc}}([0,\infty);H^{-1}_0)$.
By Definition~\ref{def:primitive}, for each $m$ we can
choose a global smooth solution $v_m$ with force
$f_{\nu_m}$ such that
\[
 \nu_m\ge m,\qquad
 \norm{v_m-u^m}_{C([0,m];H^{-1})}\le\frac1m.
\]
For fixed $T>0$ and $m\ge T$, we have
\[
 \norm{v_m-u}_{C([0,T];H^{-1})}
 \le\frac1m+\norm{u^m-u}_{C([0,T];H^{-1})}.
\]
The right-hand side tends to zero as $m\to\infty$.
Since $\nu_m\ge m$, Definition~\ref{def:primitive}
gives $u\in\cK_0^+$.

Now let $u\in\cK_0^+$ and $s\ge0$.
Choose smooth solutions $\{v_j\}$ with forces
$f_{\nu_j}$ approximating $u$ as in
Definition~\ref{def:primitive}.
Since the forcing is independent of time,
$\mathcal T_s v_j$ is again a global smooth solution
with force $f_{\nu_j}$.
For every $T>0$, convergence on $[0,s+T]$ gives
\[
 \lim_{j\to\infty}
 \norm{\mathcal T_s v_j-\mathcal T_s u}
      _{C([0,T];H^{-1})}=0.
\]
Applying Definition~\ref{def:primitive} once more,
we obtain $\mathcal T_su\in\cK_0^+$.
\end{proof}

\subsection{The trajectory attractor}

We now turn to the construction of the strong trajectory
attractor for $\cK_0^+$.
Time translations are continuous in
$C_{\mathrm{loc}}([0,\infty);H^{-1}_0)$, and
Proposition~\ref{prop:primitive-closed} ensures that they
preserve $\cK_0^+$.

The following lemma establishes the compactness of large
time translates needed for this construction.
Its first part gives a compactness criterion in $H^{-1}_0$
that will also be used for defect trajectories.

\begin{lemma}
\label{lem:asymptotic-compactness}
The following two compactness properties hold.
\leavevmode
\begin{enumerate}[
 label=\textup{(\roman*)},
 leftmargin=2.6em,
 beginpenalty=10000
]

\item
Let $Q_j\subset H^{-1}_0(\T)$ be compact for every $j\ge1$,
and assume that $\bigcup_{j\ge1}Q_j$ is bounded in $H^{-1}$.
If
\begin{equation}\label{eq:double-limit-compact-hull}
 \lim_{\kappa\to\infty}\limsup_{j\to\infty}
 \sup_{q\in Q_j}\kappa\alpha_\kappa(q)=0,
\end{equation}
then $\bigcup_{j\ge1}Q_j$ is relatively compact in $H^{-1}_0$.

\item
Let $u_j\in\cK_0^+$ and $t_j\ge0$ satisfy
\[
 \sup_{j\ge1}\norm{u_j(0)}_{H^{-1}}<\infty,
 \qquad
 \lim_{j\to\infty}t_j=\infty.
\]
Then there exist a subsequence, still indexed by $j$,
and a trajectory $v\in\cK_0^+$ such that
\[
 \lim_{j\to\infty}
 \norm{\mathcal T_{t_j}u_j-v}_{C([0,T];H^{-1})}=0
 \qquad\text{for every }T>0.
\]

\end{enumerate}
\end{lemma}

\begin{proof}
\textup{(i)} Let $Q=\bigcup_{j\ge1}Q_j$.
Since $Q$ is bounded, \eqref{eq:det-weighted-norm}
holds uniformly on $Q$ for sufficiently large $\kappa$.
For $\kappa'\ge\kappa$, we therefore have
\[
 \kappa'\alpha_{\kappa'}(q)
 \le C\norm q_{H^{-1}_{\kappa'}}^2
 \le C\norm q_{H^{-1}_{\kappa}}^2
 \le C\kappa\alpha_\kappa(q),
 \qquad q\in Q.
\]

Fix $\kappa$ sufficiently large.
For $J\ge2$ and $\kappa'\ge\kappa$, the preceding
comparison gives
\[
 \sup_{q\in Q}\kappa'\alpha_{\kappa'}(q)
 \le\max\left\{
 \sup_{q\in\bigcup_{j<J}Q_j}
       \kappa'\alpha_{\kappa'}(q),\,
 C\sup_{j\ge J}\sup_{q\in Q_j}
       \kappa\alpha_\kappa(q)
 \right\}.
\]
The finite union $\bigcup_{j<J}Q_j$ is compact.
By \Cref{thm:high-energy}, the first term inside the
maximum tends to zero as $\kappa'\to\infty$.
The second term is independent of $\kappa'$, so
\[
 0\le
 \limsup_{\kappa'\to\infty}\sup_{q\in Q}
 \kappa'\alpha_{\kappa'}(q)
 \le C\sup_{j\ge J}\sup_{q\in Q_j}
 \kappa\alpha_\kappa(q).
\]
Keeping $\kappa$ fixed and letting $J\to\infty$,
the right-hand side decreases to
$$C\limsup_{j\to\infty}\sup_{q\in Q_j}
\kappa\alpha_\kappa(q).$$
We may now let $\kappa\to\infty$ and use
\eqref{eq:double-limit-compact-hull} to obtain
\eqref{eq:beta-tight} on $Q$.
The relative compactness of $Q$ follows from
\Cref{thm:high-energy}.

\textup{(ii)}
Choose $R>0$ such that
$\sup_j\norm{u_j(0)}_{H^{-1}}\le R$, and set
$v_j=\mathcal T_{t_j}u_j$.
By Proposition~\ref{prop:primitive-closed},
$v_j\in\cK_0^+$.
Fix $T>0$. The absorbing estimate
\eqref{eq:Hminus1-absorbing} gives
\[
 \sup_j\sup_{0\le s\le T}
 \norm{v_j(s)}_{H^{-1}}<\infty.
\]
For $\kappa\ge\kappa_*(R+1,\gamma,\norm f_{H^{-1}})$,
\eqref{eq:global-beta} and
\eqref{eq:intro-spectral-comparison} yield
\begin{equation}\label{eq:shifted-beta}
 \sup_{0\le s\le T}\kappa\alpha_\kappa(v_j(s))
 \le CR^2e^{-c_0\gamma t_j}
     +C\gamma^{-2}\norm f_{H^{-1}_\kappa}^2.
\end{equation}
Taking the upper limit as $j\to\infty$, we get
\[
 0\le
 \limsup_{j\to\infty}\sup_{0\le s\le T}
 \kappa\alpha_\kappa(v_j(s))
 \le C\gamma^{-2}\norm f_{H^{-1}_\kappa}^2
 \le\frac{C}{\gamma^2\kappa^2}\norm f_{L^2}^2.
\]
Letting $\kappa\to\infty$ verifies
\eqref{eq:double-limit-compact-hull}
for the compact sets $v_j([0,T])$.
Part~\textup{(i)} therefore shows that
$\bigcup_jv_j([0,T])$ has a compact closure
$K_T$ in $H^{-1}$.

For each fixed $n\ge1$, the integrand in
\eqref{eq:primitive-Birkhoff-integral}
is continuous on $K_T$.
Integrating over $[s,t]$, we thus obtain
\[
 |\Omega_n(v_j(t))-\Omega_n(v_j(s))|
 \le C_{n,T}|t-s|,
 \qquad s,t\in[0,T],
\]
where $C_{n,T}$ is independent of $j$.
Lemma~\ref{lem:homeomorphism-equicontinuity}
now gives equicontinuity in $H^{-1}$.
The Arzel\`a--Ascoli theorem then yields relative
compactness in $C([0,T];H^{-1}_0)$.

A diagonal argument on $[0,1],[0,2],\ldots$
gives a subsequence $\{v_{j_\ell}\}$ and a path
$v\in C_{\mathrm{loc}}([0,\infty);H^{-1}_0)$ such that
\[
 \lim_{\ell\to\infty}
 \norm{v_{j_\ell}-v}_{C([0,T];H^{-1})}=0
 \qquad\text{for every }T>0.
\]
Finally, Proposition~\ref{prop:primitive-closed}
implies that $v\in\cK_0^+$.
\end{proof}

The asymptotic compactness established above allows us
to construct the strong trajectory attractor for $\cK_0^+$.

\begin{proof}[Proof of Theorem~\ref{thm:intro-primitive}]
Existence of primitive trajectories follows from
Proposition~\ref{prop:primitive-existence}.
Let $B_*\subset H^{-1}_0$ be a closed absorbing ball
as in Corollary~\ref{cor:absorbing}, and define
\begin{equation}\label{eq:primitive-omega}
 \cA=\left\{
 \lim_{j\to\infty}\mathcal T_{t_j}u_j:
 u_j\in\cK_0^+,\ u_j(0)\in B_*,
 \ \lim_{j\to\infty}t_j=\infty
 \right\},
\end{equation}
where $t_j\ge0$ and the limits are taken in
$C_{\mathrm{loc}}([0,\infty);H^{-1}_0)$.
Lemma~\ref{lem:asymptotic-compactness}\textup{(ii)}
ensures that $\varnothing\ne\cA\subset\cK_0^+$.

\emph{Compactness.}
Let $\{a_j\}\subset\cA$.
By \eqref{eq:primitive-omega}, we can choose
$u_j\in\cK_0^+$ and $t_j\ge j$ such that
\[
 u_j(0)\in B_*,
 \qquad
 \norm{a_j-\mathcal T_{t_j}u_j}_{C([0,j];H^{-1})}
 \le\frac1j.
\]
By Lemma~\ref{lem:asymptotic-compactness}\textup{(ii)},
a subsequence of $\{\mathcal T_{t_j}u_j\}$
converges to some $a\in\cA$.
The approximation above then gives, along the same subsequence,
\[
 \lim_{j\to\infty}
 \norm{a_j-a}_{C([0,T];H^{-1})}=0
 \qquad\text{for every }T>0.
\]
Thus $\cA$ is compact in the metrizable locally uniform
topology.

\emph{Attraction.}
Fix $T>0$.
If $\cA$ failed to attract trajectories starting in $B_*$,
there would be $\varepsilon>0$, times $t_j$ tending to
infinity, and trajectories $u_j\in\cK_0^+$ with
$u_j(0)\in B_*$ such that
\[
 \inf_{b\in\cA}
 \norm{\mathcal T_{t_j}u_j-b}_{C([0,T];H^{-1})}
 \ge\varepsilon.
\]
However, Lemma~\ref{lem:asymptotic-compactness}\textup{(ii)}
gives a subsequence converging to a member of $\cA$,
contradicting this inequality.

Now let $B\subset H^{-1}_0$ be bounded.
The absorbing estimate \eqref{eq:Hminus1-absorbing}
provides a common time $\tau_B\ge0$ such that
$u(\tau_B)\in B_*$, where $u\in\cK_0^+$ with $\ u(0)\in B$.
By Proposition~\ref{prop:primitive-closed},
$\mathcal T_{\tau_B}u\in\cK_0^+$.
Applying the attraction just proved to these translated
paths, with elapsed time $t-\tau_B$, gives uniform
attraction for all trajectories starting in $B$.
This is the attraction required in
Definition~\ref{def:trajectory-attractor}.

\emph{Strict invariance and uniqueness.}
Fix $s\ge0$ and $a\in\cA$, and choose a representation
as in \eqref{eq:primitive-omega}.
Continuity of time translation gives
\[
 \mathcal T_sa
 =\lim_{j\to\infty}\mathcal T_{t_j+s}u_j
 \in\cA.
\]
Hence $\mathcal T_s\cA\subset\cA$.
For the reverse inclusion, discard the finitely many
terms with $t_j<s$.
Since $t_j-s$ tends to infinity,
Lemma~\ref{lem:asymptotic-compactness}\textup{(ii)}
allows us to pass to a subsequence for which
\[
 b=\lim_{j\to\infty}\mathcal T_{t_j-s}u_j\in\cA.
\]
Another application of continuity gives
\[
 \mathcal T_sb
 =\lim_{j\to\infty}\mathcal T_{t_j}u_j=a.
\]
Thus $\mathcal T_s\cA=\cA$ for every $s\ge0$.
Uniqueness follows from strict invariance and uniform
attraction, since the initial states of any compact
trajectory attractor form a bounded set.

\emph{Complete trajectories.}
The kernel representation for trajectory attractors
\cite[Chapter~II, Theorem~4.2]{VishikChepyzhov2011},
applied to the closed trajectory class $\cK_0^+$
and its strong attractor, yields
\[
 \cA=\left\{
 U|_{[0,\infty)}:
 U\text{ is a bounded complete primitive trajectory}
 \right\}.
\]

Finally, suppose $f=0$, and let $U$ be a bounded complete
primitive trajectory.
Choose $R>0$ such that
$\sup_{s\in\mathbb R}\norm{U(s)}_{H^{-1}}\le R$.
Applying \eqref{eq:Hminus1-absorbing} to the forward
path starting at time $t-T$, we obtain
\[
 \norm{U(t)}_{H^{-1}}^2
 \le C_Re^{-c\gamma T},
 \qquad t\in\mathbb R,\quad T>0.
\]
Letting $T\to\infty$ shows that $U\equiv0$.
The kernel representation and nonemptiness of $\cA$
therefore give $\cA=\{0\}$, where $0$ denotes the
zero trajectory.
\end{proof}

\section{Defect trajectories and the global attractor}
\label{sec:defect}

We now turn from attraction in path space to the evolution
in the phase space $H^{-1}_0$.
The defect trajectories introduced in
Definition~\ref{def:defect-trajectory} are shown to generate
a strict multivalued semiflow with a compact global attractor.
We also characterize this attractor by bounded complete
defect trajectories.
In Section~\ref{sec:pullback-l2}, this characterization
will be used to prove that the attractor has $L^2$ regularity
and coincides with the classical $L^2$ attractor.

\subsection{Properties of defect trajectories}\label{5.1}

We first establish the properties of defect trajectories
needed to construct a strict multivalued semiflow.

\begin{proposition}
\label{prop:defect-structure}
The class $\cK_*^+$ has the following properties.
\begin{enumerate}[label=\textup{(\roman*)},leftmargin=2.6em]

\item
For every $u_0\in H^{-1}_0(\T)$, there exists
$u\in\cK_*^+$ such that $u(0)=u_0$.

\item
The class is preserved under forward time translations:
\[
 \mathcal T_s\cK_*^+\subseteq\cK_*^+,
 \qquad s\ge0.
\]

\item
The class is sequentially closed in
$C_{\mathrm{loc}}([0,\infty);H^{-1}_0)$.

\item
Let $T\ge0$ and $u,v\in\cK_*^+$ satisfy $u(T)=v(0)$.
Then the concatenated path
\[
 (u\diamond_Tv)(t)=
 \begin{cases}
 u(t),&0\le t\le T,\\
 v(t-T),&t>T
 \end{cases}
\]
belongs to $\cK_*^+$.

\end{enumerate}
\end{proposition}

\begin{proof}
The inclusion $\cK_0^+\subset\cK_*^+$, together with
Proposition~\ref{prop:primitive-existence}, proves
\textup{(i)}.

For \textup{(ii)}, let $u\in\cK_*^+$ and fix $s\ge0$.
By Definition~\ref{def:defect-trajectory}, we can choose
piecewise classical approximants $w_j$ on $[0,s+j]$
such that
\[
 \norm{w_j-u}_{C_*([0,s+j];H^{-1})}\le j^{-1},
 \qquad
 \mathfrak d_{s+j}(w_j)\le j^{-1},
 \qquad
 \nu_*(w_j)\ge j.
\]
Set
\[
 \widetilde w_j(t)=w_j(s+t),
 \qquad 0\le t\le j,
\]
using the right trace at the new initial time.
Restriction introduces no additional jumps and can only
increase the minimum force cutoff.
Thus the same three bounds hold on $[0,j]$ for
$\widetilde w_j$ relative to $\mathcal T_su$.
Definition~\ref{def:defect-trajectory} therefore gives
$\mathcal T_su\in\cK_*^+$.

For \textup{(iii)}, let $\{u^m\}\subset\cK_*^+$ converge
to $u$ in $C_{\mathrm{loc}}([0,\infty);H^{-1}_0)$.
For each $m$, choose a piecewise classical approximant
$w_m$ on $[0,m]$ satisfying
\begin{equation}\label{eq:closedness-diagonal-choice}
 \norm{w_m-u^m}_{C_*([0,m];H^{-1})}\le m^{-1},
 \quad
 \mathfrak d_m(w_m)\le m^{-1},
 \quad
 \nu_*(w_m)\ge m.
\end{equation}
For fixed $T>0$ and $m\ge T$, we have
\[
 \norm{w_m-u}_{C_*([0,T];H^{-1})}
 \le m^{-1}+\norm{u^m-u}_{C([0,T];H^{-1})}.
\]
The right-hand side tends to zero as $m\to\infty$.
Together with \eqref{eq:closedness-diagonal-choice},
this verifies all three conditions in
\eqref{eq:defect-convergence}, so $u\in\cK_*^+$.

For \textup{(iv)}, the case $T=0$ is immediate.
Suppose $T>0$, and choose approximants $a_j$ on $[0,T]$
and $b_j$ on $[0,j]$ such that
\[
 \begin{gathered}
 \norm{a_j-u}_{C_*([0,T];H^{-1})}\le j^{-1},
 \qquad
 \norm{b_j-v}_{C_*([0,j];H^{-1})}\le j^{-1},\\
 \mathfrak d_T(a_j)\le j^{-1},
 \qquad
 \mathfrak d_j(b_j)\le j^{-1},
 \qquad
 \min\{\nu_*(a_j),\nu_*(b_j)\}\ge j.
 \end{gathered}
\]
At time $T$, join $a_j$ to $t\mapsto b_j(t-T)$.
Since $u(T)=v(0)$, the additional jump satisfies
\[
 \begin{aligned}
 \norm{b_j(0+)-a_j(T-)}_{H^{-1}}
 \le\norm{b_j(0+)-v(0)}_{H^{-1}} +\norm{u(T)-a_j(T-)}_{H^{-1}}
 \le2j^{-1}.
 \end{aligned}
\]
Let $c_j$ be the resulting path restricted to $[0,j]$.
Adding the two restart errors and the new jump, we obtain
\[
 \mathfrak d_j(c_j)\le4j^{-1},
 \qquad
 \nu_*(c_j)\ge j.
\]
Moreover, for every fixed $L>0$ and all sufficiently
large $j$,
\[
 \norm{c_j-u\diamond_Tv}_{C_*([0,L];H^{-1})}
 \le j^{-1}.
\]
These bounds verify \eqref{eq:defect-convergence}
for $u\diamond_Tv$, proving \textup{(iv)}.
\end{proof}

The number of restarts need not be uniformly bounded along
the approximating sequence.
We therefore control their cumulative effect on
$\alpha_\kappa$ by the total $H^{-1}$ jump size
$\mathfrak d_T(w)$.
The following lemma gives the estimate for a single restart.

\begin{lemma}
\label{lem:defect-alpha}
Fix $\kappa\ge1$, $R>0$, and $0<\eta<1$.
Let $q,J\in H^{-1}_0(\T)$ satisfy
\[
 \norm{q+\theta J}_{H^{-1}}\le R,
 \qquad
 \norm{A_\kappa(q+\theta J)}_{\op}\le\eta
 \qquad\text{for every }\theta\in[0,1].
\]
Then
\begin{equation}\label{eq:alpha-jump}
 \abs{\alpha_\kappa(q+J)-\alpha_\kappa(q)}
 \le C_{\kappa,\eta,R}\norm J_{H^{-1}},
\end{equation}
where the constant depends only on $\kappa$, $\eta$, and $R$.
\end{lemma}

\begin{proof}
Integrating the derivative bound
\eqref{eq:alpha-local-Lipschitz} along $q+\theta J$,
we obtain
\[
 \begin{aligned}
 \abs{\alpha_\kappa(q+J)-\alpha_\kappa(q)}
 &=\left|\int_0^1
 D\alpha_\kappa(q+\theta J)[J]\,\dd\theta\right|\\
 &\le C_{\kappa,\eta,R}\norm J_{H^{-1}}.
 \end{aligned}
\]
This completes the proof.
\end{proof}

We next extend the dissipative estimate to paths with
finitely many restarts.
The bounds will depend on the total jump size and remain
uniform in the number of pieces.

Fix $0<\eta<1$, $T>0$, and $\kappa\ge1$.
Let $w$ be a piecewise smooth path on the partition
\[
 0=t_0<\cdots<t_N=T,
\]
and write
\[
 J_\ell=w(t_\ell^+)-w(t_\ell^-),
 \qquad
 \mathfrak d_T(w)
 =\sum_{\ell=1}^{N-1}\norm{J_\ell}_{H^{-1}}.
\]
At internal restart times we use the right trace,
$w(t_\ell)=w(t_\ell^+)$, and set $w(T)=w(T^-)$.
Assume that all states and connecting segments lie in
a fixed $H^{-1}$ ball, and put
\[
 C_\kappa
 =\sup_{\norm q_{H^{-1}}\le1}
   \norm{A_\kappa(q)}_{\op}.
\]

Suppose that, on each smooth piece, whenever
$\norm{A_\kappa(w(t))}_{\op}<\eta$, we have
\begin{equation}\label{eq:restart-differential}
 \frac{\dd}{\dd t}\alpha_\kappa(w(t))
 +\lambda\alpha_\kappa(w(t))\le F_0
\end{equation}
and
\begin{equation}\label{eq:restart-coercivity}
 \alpha_\kappa(w(t))
 \ge c_\eta\norm{A_\kappa(w(t))}_{\HS}^2.
\end{equation}
Here $\lambda>0$, $F_0\in\mathbb R$, and $c_\eta>0$
are common to all pieces.

At a restart, write $q=w(t_\ell^-)$ and $J=J_\ell$.
Whenever
\[
 \norm{A_\kappa(q+\theta J)}_{\op}\le\frac{3\eta}{4}
 \qquad\text{for every }\theta\in[0,1],
\]
Lemma~\ref{lem:defect-alpha} gives
\begin{equation}\label{eq:restart-jump-Lipschitz}
 \abs{\alpha_\kappa(q+J)-\alpha_\kappa(q)}
 \le L\norm J_{H^{-1}},
\end{equation}
with a constant $L\ge0$ common to all restarts.

Assume the initial and source bounds
\begin{equation}\label{eq:restart-initial-margins}
 \norm{A_\kappa(w(0))}_{\op}\le\frac{\eta}{4},
 \qquad
 \max\left\{
 \alpha_\kappa(w(0)),\frac{F_0}{\lambda}
 \right\}
 \le\frac{c_\eta\eta^2}{8},
\end{equation}
together with
\begin{equation}\label{eq:restart-total-defect-margins}
 C_\kappa\mathfrak d_T(w)\le\frac{\eta}{4},
 \qquad
 L\mathfrak d_T(w)\le\frac{c_\eta\eta^2}{8}.
\end{equation}
We now check that the estimates remain applicable
throughout the path.

Starting from $w(0)$, we proceed along the smooth pieces
and across the restarts.
On the portion already reached by the argument,
integration and the jump estimate give
\[
 \alpha_\kappa(w(t))
 \le
 \max\left\{
 \alpha_\kappa(w(0)),\frac{F_0}{\lambda}
 \right\}
 +L\sum_{0<t_\ell\le t}\norm{J_\ell}_{H^{-1}}
 \le\frac{c_\eta\eta^2}{4}.
\]
By \eqref{eq:restart-coercivity},
$\norm{A_\kappa(w(t))}_{\op}\le\eta/2$.
Continuity therefore prevents the operator norm from
reaching $\eta$ on the current smooth piece.

At the next restart, the connecting segment satisfies
\[
 \begin{aligned}
 \norm{A_\kappa(q+\theta J)}_{\op}
 \le\norm{A_\kappa(q)}_{\op}
      +C_\kappa\norm J_{H^{-1}}\le\frac{\eta}{2}+\frac{\eta}{4}
 =\frac{3\eta}{4},
 \qquad 0\le\theta\le1.
 \end{aligned}
\]
Thus \eqref{eq:restart-jump-Lipschitz} applies, and adding
the new jump to the cumulative sum preserves the same
energy bound.
Since $\norm{A_\kappa(w(t_\ell^+))}_{\op}<\eta$,
\eqref{eq:restart-coercivity} holds at the right trace
by continuity from the next smooth piece.
The operator norm is consequently bounded by $\eta/2$
there as well.

Repeating this argument gives
\begin{equation}\label{eq:restart-invariant-region}
 \norm{A_\kappa(w(t^\pm))}_{\op}\le\frac{\eta}{2},
 \qquad
 \alpha_\kappa(w(t^\pm))\le\frac{c_\eta\eta^2}{4}
\end{equation}
at every ordinary time and every available one-sided trace.
In particular, the determinant is defined throughout
$[0,T]$, and the differential and jump estimates can
be used on any subinterval.

\begin{auxiliarylemma}
\label{lem:arbitrary-restart-induction}
Under the assumptions above, for every $0\le s\le t\le T$, we have
\begin{equation}\label{eq:restart-Gronwall}
 \alpha_\kappa(w(t))
 \le e^{-\lambda(t-s)}\alpha_\kappa(w(s))
 +\frac{F_0}{\lambda}
  \bigl(1-e^{-\lambda(t-s)}\bigr)
 +L\mathfrak d_T(w).
\end{equation}
\end{auxiliarylemma}

\begin{proof}
Integrating \eqref{eq:restart-differential} on each smooth
piece and applying \eqref{eq:restart-jump-Lipschitz}
at the intervening restarts, we obtain
\[
 \begin{aligned}
 \alpha_\kappa(w(t))
 \le e^{-\lambda(t-s)}\alpha_\kappa(w(s))
 +\frac{F_0}{\lambda}
  \bigl(1-e^{-\lambda(t-s)}\bigr)
 +L\sum_{s<t_\ell\le t}
 e^{-\lambda(t-t_\ell)}\norm{J_\ell}_{H^{-1}}.
 \end{aligned}
\]
Since $\lambda>0$, the weighted sum is at most
$\mathfrak d_T(w)$.
This completes the proof.
\end{proof}

We can now extend the determinant estimates and the
absorbing bound to defect trajectories.

\begin{proposition}
\label{prop:defect-beta}
Let $B\subset H^{-1}_0$ be bounded.
There exists $\kappa_B\ge1$, depending only on an
$H^{-1}$ bound for $B$, $\gamma$, and $\norm f_{H^{-1}}$,
with the following property.
For every $u\in\cK_*^+$ with $u(0)\in B$ and every
$\kappa\ge\kappa_B$, the determinant is defined along $u$,
and \eqref{eq:intro-spectral-comparison} holds.
Moreover, for every $0\le s\le t$,
\begin{equation}\label{eq:defect-limit-beta}
 \kappa\alpha_\kappa(u(t))
 \le e^{-c\gamma(t-s)}\kappa\alpha_\kappa(u(s))
 +C\gamma^{-2}\norm f_{H^{-1}_\kappa}^2
  \bigl(1-e^{-c\gamma(t-s)}\bigr).
\end{equation}
Here $c,C>0$ are absolute constants.
The absorbing estimate \eqref{eq:Hminus1-absorbing}
also holds uniformly over these trajectories.
\end{proposition}

\begin{proof}
Choose $R>0$ such that
$\sup_{q\in B}\norm q_{H^{-1}}\le R$, and fix $\eta=1/4$.
By \eqref{eq:alpha-coercive} and
\eqref{eq:HS-weighted-norm}, the smooth pieces satisfy
\eqref{eq:restart-differential} with
\[
 \lambda=c\gamma,
 \qquad
 F_0=C\gamma^{-1}\kappa^{-1}
       \norm f_{H^{-1}_\kappa}^2,
\]
uniformly over all Fourier cutoffs of $f$.
The lower bound in \eqref{eq:alpha-HS-comparison}
gives \eqref{eq:restart-coercivity}.
The upper bound, together with
\eqref{eq:HS-weighted-norm}, allows us to choose
$\kappa_B$ so that \eqref{eq:restart-initial-margins}
holds whenever
\[
 \norm q_{H^{-1}}\le R+1,
 \qquad \kappa\ge\kappa_B.
\]
This choice depends only on $R$, $\gamma$, and
$\norm f_{H^{-1}}$.

Now fix $u\in\cK_*^+$ with $u(0)\in B$, $T>0$, and
$\kappa\ge\kappa_B$.
Choose piecewise classical approximants $w_j$ as in
Definition~\ref{def:defect-trajectory}, restricted to $[0,T]$.
Then
\[
 \lim_{j\to\infty}
 \norm{w_j-u}_{C_*([0,T];H^{-1})}=0,
 \qquad
 \lim_{j\to\infty}\mathfrak d_T(w_j)=0.
\]
After discarding finitely many terms,
$\norm{w_j(0)}_{H^{-1}}\le R+1$.
The $C_*$ convergence also places all states, traces,
and connecting segments in one fixed $H^{-1}$ ball.
Lemma~\ref{lem:defect-alpha} therefore supplies
\eqref{eq:restart-jump-Lipschitz} with a common constant
$L_{\kappa,T,u}$.

At this fixed $\kappa$, both $C_\kappa$ and
$L_{\kappa,T,u}$ are independent of $j$.
Thus \eqref{eq:restart-total-defect-margins} holds for
all sufficiently large $j$.
Consequently, \eqref{eq:restart-invariant-region} holds
for $w_j$, and Lemma~\ref{lem:arbitrary-restart-induction}
gives, after multiplication by $\kappa$,
\begin{align}
 \kappa\alpha_\kappa(w_j(t))
 &\le e^{-c\gamma(t-s)}\kappa\alpha_\kappa(w_j(s))
 +C\gamma^{-2}\norm f_{H^{-1}_\kappa}^2
  \bigl(1-e^{-c\gamma(t-s)}\bigr)\notag\\
 &\quad+\kappa L_{\kappa,T,u}\mathfrak d_T(w_j),
 \qquad 0\le s\le t\le T.
 \label{eq:defect-prelimit-beta}
\end{align}

By \eqref{eq:HS-weighted-norm}, the operator bound
passes to $u$.
The determinant is therefore defined along $u$, and
\eqref{eq:det-weighted-norm} gives
\eqref{eq:intro-spectral-comparison}.
Applying Lemma~\ref{lem:defect-alpha} between
$w_j(r)$ and $u(r)$, we also obtain
\[
 \lim_{j\to\infty}\sup_{0\le r\le T}
 \abs{\alpha_\kappa(w_j(r))-\alpha_\kappa(u(r))}=0.
\]
Since $\kappa L_{\kappa,T,u}\mathfrak d_T(w_j)$
tends to zero at this fixed $\kappa$, passing to the
limit in \eqref{eq:defect-prelimit-beta} proves
\eqref{eq:defect-limit-beta}.

As $T>0$ and $\kappa\ge\kappa_B$ were arbitrary,
the estimate holds on $[0,\infty)$ for every such $\kappa$.
The approximation index may depend on $\kappa$, $T$,
and $u$, while $\kappa_B$ and the limiting constants
are uniform for $u(0)\in B$.
Finally, the fixed-$\kappa$ argument in
Corollary~\ref{cor:absorbing} gives
\eqref{eq:Hminus1-absorbing} with the same uniformity.
This completes the proof.
\end{proof}

We next establish the integral Birkhoff coordinate equations
for defect trajectories.
These identities will be used to prove time equicontinuity
below and to identify the $L^2$ dynamics in
Section~\ref{sec:pullback-l2}.

\begin{theorem}
\label{thm:birkhoff-integral}
Let $u\in\cK_*^+$.
Then, for every $n\ge1$ and $0\le s\le t$,
\begin{equation}\label{eq:birkhoff-integral}
 \Omega_n(u(t))-\Omega_n(u(s))
 =\int_s^t
 \left[
 i\omega_n(u(\tau))\Omega_n(u(\tau))
 +D\Omega_n(u(\tau))[f-\gamma u(\tau)]
 \right]\dd\tau.
\end{equation}
\end{theorem}

\begin{proof}
Fix $n\ge1$ and $0\le s\le t<T$.
Choose piecewise classical approximants $w_j$ for $u$
as in Definition~\ref{def:defect-trajectory},
restricted to $[0,T]$, and denote their piecewise
constant forcing by $g_j$.
On each smooth piece, the chain rule and
\eqref{eq:Birkhoff-KdV-vector} give
\begin{equation}\label{eq:piece-coordinate-ODE}
 \frac{\dd}{\dd\tau}\Omega_n(w_j)
 =i\omega_n(w_j)\Omega_n(w_j)
  +D\Omega_n(w_j)[g_j-\gamma w_j].
\end{equation}
Using the right trace at internal restart times,
we integrate over the pieces intersecting $[s,t]$
and obtain
\begin{align}
 \Omega_n(w_j(t))-\Omega_n(w_j(s))
 &=\int_s^t
 \left[
 i\omega_n(w_j)\Omega_n(w_j)
 +D\Omega_n(w_j)[g_j-\gamma w_j]
 \right]\dd\tau\notag\\
 &\quad+\sum_{t_\ell\in(s,t]}
 \left[
 \Omega_n(w_j(t_\ell^+))
 -\Omega_n(w_j(t_\ell^-))
 \right].
 \label{eq:piece-summed-coordinate}
\end{align}

Let $K_T$ be the closure in $H^{-1}$ of $u([0,T])$
together with all states, one-sided traces, and joining
segments of the paths $w_j$.
Each approximant contributes finitely many compact
piece images and segments.
Moreover, the $C_*$ convergence places all these sets
uniformly close to $u([0,T])$ as $j\to\infty$.
Thus $K_T$ is compact.

Continuity of $D\Omega_n$ on $K_T$ gives a common
operator norm bound.
Integrating along each joining segment, we get
\begin{equation}\label{eq:Birkhoff-jumps-vanish}
 \sum_{t_\ell\in(s,t]}
 \abs{\Omega_n(w_j(t_\ell^+))
      -\Omega_n(w_j(t_\ell^-))}
 \le C_{n,T}\mathfrak d_T(w_j),
\end{equation}
where $C_{n,T}$ is independent of $j$.
Since $\mathfrak d_T(w_j)$ tends to zero,
the jump sum vanishes in the limit.

For the forcing, we have
\begin{equation}\label{eq:piece-force-uniform}
 \sup_{0\le\tau\le T}\norm{g_j(\tau)-f}_{H^{-1}}
 \le\sup_{\nu\ge\nu_*(w_j)}
       \norm{P_{\le\nu}f-f}_{H^{-1}}.
\end{equation}
The right-hand side tends to zero because
$\nu_*(w_j)$ tends to infinity.

Finally, $\Omega_n$, $\omega_n$, and $D\Omega_n$
are uniformly continuous on $K_T$, with operator norm
used for $D\Omega_n$.
Together with the path and forcing convergence,
this gives uniform convergence of the coordinate
values and integrands in
\eqref{eq:piece-summed-coordinate}.
Passing to the limit and using
\eqref{eq:Birkhoff-jumps-vanish}, we obtain
\eqref{eq:birkhoff-integral}.
This completes the proof.
\end{proof}
 
We next show that the smooth approximating equations
determine a distributional limit of the quadratic term.

Fix $u\in\cK_*^+$ and $T>0$.
Let $w_j$ be piecewise classical approximants for $u$
as in Definition~\ref{def:defect-trajectory},
restricted to $[0,T]$, with partitions
$
 0=t_{j,0}<\cdots<t_{j,N_j}=T.
$
Denote their piecewise constant forces by $g_j$,
using the same endpoint convention as for $w_j$.
Then
\begin{equation}\label{eq:defect-force-convergence-distribution}
 \lim_{j\to\infty}
 \sup_{0\le t\le T}\norm{g_j(t)-f}_{H^{-1}}=0,
 \qquad
 \lim_{j\to\infty}\mathfrak d_T(w_j)=0.
\end{equation}

The jumps contribute the $H^{-1}$-valued measure
\[
 \mathcal J_j
 =\sum_{\ell=1}^{N_j-1}
 \bigl(w_j(t_{j,\ell}^+)-w_j(t_{j,\ell}^-)\bigr)
 \delta_{t_{j,\ell}},
\]
where $\delta_{t_{j,\ell}}$ denotes the Dirac mass in time.
The distributional time derivative contains this measure,
so the classical equations on the pieces give
\begin{equation}\label{eq:piecewise-equation-with-jumps}
 \partial_tw_j+\partial_x^3w_j
 -3\partial_x(w_j^2)+\gamma w_j
 =g_j+\mathcal J_j
 \quad\text{in }\mathcal D'((0,T)\times\T).
\end{equation}
For every $\varphi\in C_c^\infty((0,T)\times\T)$,
Sobolev duality yields
\begin{equation}\label{eq:jump-measure-vanishes}
 \abs{\ip{\mathcal J_j}{\varphi}}
 \le 2\pi\,\mathfrak d_T(w_j)
 \sup_{0\le t\le T}
 \norm{\varphi(t,\cdot)}_{H^1}.
\end{equation}
Thus $\mathcal J_j$ converges to zero in
$\mathcal D'((0,T)\times\T)$.

For the limiting path, define
\begin{equation}\label{eq:pathwise-defect-nonlinearity}
 \mathfrak N[u]
 :=\frac13\bigl(
 \partial_tu+\partial_x^3u+\gamma u-f
 \bigr).
\end{equation}
All terms on the right are well defined as distributions,
so $\mathfrak N[u]\in
\mathcal D'((0,\infty)\times\T)$.

\begin{theorem}
\label{thm:closed-nonlinear-graph}
Let $u\in\cK_*^+$.
For every $T>0$, any sequence $w_j$ approximating $u$
as in Definition~\ref{def:defect-trajectory} satisfies
\begin{equation}\label{eq:closed-product-limit}
 \lim_{j\to\infty}\partial_x(w_j^2)
 =\mathfrak N[u]
 \quad\text{in }\mathcal D'((0,T)\times\T),
\end{equation}
where $\mathfrak N[u]$ is defined by
\eqref{eq:pathwise-defect-nonlinearity}.
\end{theorem}

\begin{proof}
Rearranging \eqref{eq:piecewise-equation-with-jumps},
we have
\[
 \partial_x(w_j^2)
 =\frac13\bigl(
 \partial_tw_j+\partial_x^3w_j
 +\gamma w_j-g_j-\mathcal J_j
 \bigr).
\]
The path convergence, together with
\eqref{eq:defect-force-convergence-distribution}
and \eqref{eq:jump-measure-vanishes}, allows us to
pass to the limit in every term on the right.
The resulting limit is $\mathfrak N[u]$ by
\eqref{eq:pathwise-defect-nonlinearity}.
This completes the proof.
\end{proof}

Consequently, every defect trajectory satisfies
\begin{equation}\label{eq:closed-graph-KdV}
 \partial_tu+\partial_x^3u
 -3\mathfrak N[u]+\gamma u=f
 \quad\text{in }\mathcal D'((0,\infty)\times\T).
\end{equation}

The argument also applies to other approximations
whose pieces solve the classical forced equation.
On every finite time interval, it requires convergence
to the same path $u$ in $C_*$, uniform convergence
of the forces to $f$ in $H^{-1}$, and vanishing total
restart error.
Thus $\mathfrak N[u]$ is independent of the partition
and the chosen approximants.
This representation independence concerns a fixed path
and does not imply uniqueness of the evolution.

Identifying $\mathfrak N[u]$ with $\partial_x(u^2)$
requires a separate argument, as illustrated by
\eqref{eq:quadratic-discontinuity}.
For forward defect trajectories with $L^2$ initial data
and bounded complete defect trajectories, this
identification is established in
Corollary~\ref{cor:complete-standard-distributional}.

\subsection{Multivalued semiflow and attractor}\label{5.2}

We now show that the evolution defined in
\eqref{eq:multiflow} is a strict multivalued semiflow
on $H^{-1}_0$, with a closed graph at each fixed time.

\begin{theorem}\label{thm:strict-multiflow}
For every $u_0\in H^{-1}_0(\T)$ and $s,t\ge0$,  we have
\[
 \cS_*(t)u_0\ne\varnothing,
 \qquad
 \cS_*(0)u_0=\{u_0\},
\]
and
\[
 \cS_*(t+s)u_0
 =\cS_*(t)\bigl(\cS_*(s)u_0\bigr)
 =\bigcup_{y\in\cS_*(s)u_0}\cS_*(t)y.
\]
Thus $(\cS_*(t))_{t\ge0}$ is a strict multivalued
semiflow on $H^{-1}_0(\T)$.
Moreover, for every fixed $t\ge0$, its graph is closed
in the strong $H^{-1}\times H^{-1}$ topology.
\end{theorem}

\begin{proof}
\emph{Strict semiflow property.}
Nonemptiness and the identity at zero follow from
Proposition~\ref{prop:defect-structure}\textup{(i)}
and \eqref{eq:multiflow}.

Let $z\in\cS_*(t+s)u_0$.
Choose $u\in\cK_*^+$ with $u(0)=u_0$ and $u(t+s)=z$.
By Proposition~\ref{prop:defect-structure}\textup{(ii)},
$\mathcal T_su\in\cK_*^+$, so
\[
 u(s)\in\cS_*(s)u_0,
 \qquad
 z\in\cS_*(t)u(s).
\]
This proves one inclusion in the composition law.

Conversely, let $y\in\cS_*(s)u_0$ and
$z\in\cS_*(t)y$.
There are $u,v\in\cK_*^+$ such that
\[
 u(0)=u_0,\qquad u(s)=v(0)=y,\qquad v(t)=z.
\]
Proposition~\ref{prop:defect-structure}\textup{(iv)}
gives $u\diamond_s v\in\cK_*^+$.
This path starts at $u_0$ and reaches $z$ at time $s+t$,
hence $z\in\cS_*(s+t)u_0$.
The reverse inclusion follows, proving strictness.

\emph{Closed graph.}
Fix $t\ge0$.
Suppose that $y_j\in\cS_*(t)u_{0,j}$ and
\[
 \lim_{j\to\infty}\norm{u_{0,j}-u_0}_{H^{-1}}=0,
 \qquad
 \lim_{j\to\infty}\norm{y_j-y}_{H^{-1}}=0.
\]
Choose $u_j\in\cK_*^+$ with
$u_j(0)=u_{0,j}$ and $u_j(t)=y_j$.
The set
\[
 Q=\{u_0\}\cup\{u_{0,j}:j\ge1\}
\]
is compact in $H^{-1}_0$.
By Proposition~\ref{prop:defect-beta},
\[
 \sup_j\sup_{\sigma\ge0}
 \norm{u_j(\sigma)}_{H^{-1}}<\infty.
\]

Fix an integer $M\ge1$.
For $\kappa\ge\kappa_Q$, where $\kappa_Q$ is given by
Proposition~\ref{prop:defect-beta} with $B=Q$,
\eqref{eq:defect-limit-beta} yields
\[
 0\le
 \sup_j\sup_{0\le\sigma\le M}
 \kappa\alpha_\kappa(u_j(\sigma))
 \le
 \sup_{q\in Q}\kappa\alpha_\kappa(q)
 +\frac{C}{\gamma^2\kappa^2}\norm f_{L^2}^2.
\]
Since $Q$ is compact, the right-hand side tends to zero
as $\kappa\to\infty$ by \Cref{thm:high-energy}.
Applying the same criterion to the bounded state set
\[
 \{u_j(\sigma):j\ge1,\ 0\le\sigma\le M\},
\]
we find that its closure $K_M$ in $H^{-1}$ is compact.

For each fixed $n\ge1$, the integrand in
\eqref{eq:birkhoff-integral} is bounded on $K_M$.
Thus
\[
 |\Omega_n(u_j(s))-\Omega_n(u_j(r))|
 \le C_{n,M}|s-r|,
 \qquad r,s\in[0,M],
\]
with $C_{n,M}$ independent of $j$.
Lemma~\ref{lem:homeomorphism-equicontinuity}
gives equicontinuity in $H^{-1}$, and the
Arzel\`a--Ascoli theorem  yields relative compactness in
$C([0,M];H^{-1}_0)$.

A diagonal argument on $[0,1],[0,2],\ldots$
now gives a subsequence $\{u_{j_\ell}\}$ and
$u\in C_{\mathrm{loc}}([0,\infty);H^{-1}_0)$ such that
\[
 \lim_{\ell\to\infty}
 \norm{u_{j_\ell}-u}_{C([0,M];H^{-1})}=0
 \qquad\text{for every integer }M\ge1.
\]
Proposition~\ref{prop:defect-structure}\textup{(iii)}
then implies $u\in\cK_*^+$.
Passing to the limit at times $0$ and $t$, we obtain
\[
 u(0)=\lim_{\ell\to\infty}u_{0,j_\ell}=u_0,
 \qquad
 u(t)=\lim_{\ell\to\infty}y_{j_\ell}=y.
\]
Hence $y\in\cS_*(t)u_0$, proving that the graph is closed.
This completes the proof.
\end{proof}
To construct the global attractor, we next establish
the existence of a bounded absorbing set and asymptotic
compactness.

\begin{proposition}\label{prop:multiflow-asymptotic}
The semiflow $(\cS_*(t))_{t\ge0}$ admits a bounded
absorbing set in $H^{-1}_0$ and is asymptotically compact.
More precisely, let $u_{0,j}\in H^{-1}_0$, $t_j\ge0$,
and $y_j\in\cS_*(t_j)u_{0,j}$ satisfy
\[
 \sup_j\norm{u_{0,j}}_{H^{-1}}<\infty,
 \qquad
 \lim_{j\to\infty}t_j=\infty.
\]
Then there exist a subsequence $\{y_{j_\ell}\}$ and
$y\in H^{-1}_0$ such that
\[
 \lim_{\ell\to\infty}\norm{y_{j_\ell}-y}_{H^{-1}}=0.
\]
\end{proposition}

\begin{proof}
The bounded absorbing set follows from
Proposition~\ref{prop:defect-beta}.
Choose $R>0$ such that
$\sup_j\norm{u_{0,j}}_{H^{-1}}\le R$, and take
$u_j\in\cK_*^+$ with
\[
 u_j(0)=u_{0,j},\qquad u_j(t_j)=y_j.
\]
The absorbing estimate also gives
\[
 \sup_j\norm{y_j}_{H^{-1}}<\infty.
\]

For all sufficiently large $\kappa$,
\eqref{eq:defect-limit-beta} and
\eqref{eq:intro-spectral-comparison} yield
\begin{equation}\label{eq:endpoint-asymptotic-beta}
 \kappa\alpha_\kappa(y_j)
 \le CR^2e^{-c\gamma t_j}
     +C\gamma^{-2}\norm f_{H^{-1}_\kappa}^2.
\end{equation}
Keeping $\kappa$ fixed and taking the upper limit
as $j\to\infty$, we obtain
\[
 0\le\limsup_{j\to\infty}\kappa\alpha_\kappa(y_j)
 \le C\gamma^{-2}\norm f_{H^{-1}_\kappa}^2
 \le\frac{C}{\gamma^2\kappa^2}\norm f_{L^2}^2.
\]
Letting $\kappa\to\infty$ verifies
\eqref{eq:double-limit-compact-hull}
for the compact singleton sets $\{y_j\}$.
Lemma~\ref{lem:asymptotic-compactness}\textup{(i)}
therefore shows that $\{y_j:j\ge1\}$ is relatively
compact in $H^{-1}_0$.
Hence $\{y_j\}$ has a strongly convergent subsequence.
This completes the proof.
\end{proof}

We now construct the global attractor and describe it
in terms of bounded complete defect trajectories.
Let $B_*\subset H^{-1}_0$ be a closed bounded absorbing
ball, and define
\begin{equation}\label{eq:multivalued-attractor-definition}
 \cA_*
 =\bigcap_{\tau\ge0}
  \overline{\bigcup_{t\ge\tau}\cS_*(t)B_*},
\end{equation}
where the closure is taken in $H^{-1}$.

\begin{theorem}
\label{thm:multivalued-attractor}
The set $\cA_*$ defined in
\eqref{eq:multivalued-attractor-definition}
is nonempty and compact in $H^{-1}_0$.
It attracts every bounded subset of the phase space:
\begin{equation}\label{eq:state-attraction}
 \lim_{t\to\infty}
 \dist_{H^{-1}}(\cS_*(t)B,\cA_*)=0
 \qquad\text{for every bounded }B\subset H^{-1}_0.
\end{equation}
Moreover,
\begin{equation}\label{eq:multivalued-strict-invariance}
 \cS_*(t)\cA_*=\cA_*,
 \qquad t\ge0,
\end{equation}
and
\begin{equation}\label{eq:defect-kernel-representation}
 \cA_*=\left\{
 U(0):U\text{ is a bounded complete defect trajectory}
 \right\}.
\end{equation}
Thus $\cA_*$ is a compact global attractor for
$(\cS_*(t))_{t\ge0}$.
\end{theorem}

\begin{proof}
The definition gives the sequential characterization
\begin{equation}\label{eq:multivalued-omega-sequential}
 \cA_*=\left\{
 \lim_{j\to\infty}y_j:
 y_j\in\cS_*(t_j)B_*,
 \ \lim_{j\to\infty}t_j=\infty
 \right\},
\end{equation}
with limits taken in $H^{-1}$.

\emph{Compactness and attraction.}
Choose $b\in B_*$ and $y_j\in\cS_*(j)b$.
Proposition~\ref{prop:multiflow-asymptotic}
gives a convergent subsequence, whose limit belongs
to $\cA_*$ by \eqref{eq:multivalued-omega-sequential}.
Hence $\cA_*\ne\varnothing$.

Now let $\{a_j\}\subset\cA_*$.
Using \eqref{eq:multivalued-omega-sequential}, choose
\[
 t_j\ge j,\qquad
 y_j\in\cS_*(t_j)B_*,
 \qquad
 \norm{y_j-a_j}_{H^{-1}}\le j^{-1}.
\]
Asymptotic compactness gives a subsequence of $\{y_j\}$
converging to some $a\in\cA_*$.
Along the same subsequence, $a_j$ has the same limit.
Thus $\cA_*$ is compact.

If $\cA_*$ failed to attract $B_*$, there would be
$\varepsilon>0$, times $t_j$ tending to infinity,
and $y_j\in\cS_*(t_j)B_*$ such that
\[
 \dist_{H^{-1}}(\{y_j\},\cA_*)\ge\varepsilon.
\]
A subsequence converges to a member of $\cA_*$,
contradicting this inequality.

For any bounded $B\subset H^{-1}_0$, choose $T_B\ge0$
such that $\cS_*(T_B)B\subseteq B_*$.
Strictness gives
\[
 \cS_*(t)B
 =\cS_*(t-T_B)\bigl(\cS_*(T_B)B\bigr)
 \subseteq\cS_*(t-T_B)B_*,
 \qquad t\ge T_B.
\]
Attraction of $B_*$ therefore yields
\eqref{eq:state-attraction}.

\emph{Strict invariance.}
Fix $t\ge0$ and $a\in\cA_*$.
Choose $y_j\in\cS_*(t_j)B_*$ converging to $a$,
with $t_j$ tending to infinity.
For all sufficiently large $j$, strictness allows us
to choose
\[
 b_j\in\cS_*(t_j-t)B_*,
 \qquad
 y_j\in\cS_*(t)b_j.
\]
By Proposition~\ref{prop:multiflow-asymptotic},
a subsequence of $\{b_j\}$ converges to some
$b\in\cA_*$.
Since $\cS_*(t)$ has a closed graph,
\[
 a\in\cS_*(t)b\subseteq\cS_*(t)\cA_*.
\]
Thus $\cA_*\subseteq\cS_*(t)\cA_*$ for every $t\ge0$.

For the reverse inclusion, this result and strictness
give, for every $\tau\ge0$,
\[
 \cS_*(t)\cA_*
 \subseteq
 \cS_*(t)\bigl(\cS_*(\tau)\cA_*\bigr)
 =\cS_*(t+\tau)\cA_*.
\]
Consequently,
\[
 \dist_{H^{-1}}(\cS_*(t)\cA_*,\cA_*)
 \le
 \dist_{H^{-1}}(\cS_*(t+\tau)\cA_*,\cA_*).
\]
Since $\cA_*$ is bounded, the right-hand side tends
to zero as $\tau\to\infty$ by
\eqref{eq:state-attraction}.
Closedness of $\cA_*$ now gives
$\cS_*(t)\cA_*\subseteq\cA_*$,
proving \eqref{eq:multivalued-strict-invariance}.

\emph{Kernel representation.}
Fix $a\in\cA_*$ and choose $v_0\in\cK_*^+$
with $v_0(0)=a$.
Using $\cS_*(1)\cA_*=\cA_*$, we can recursively
choose $v_m\in\cK_*^+$, $m\ge1$, such that
\[
 v_m(0)\in\cA_*,
 \qquad
 v_m(1)=v_{m-1}(0).
\]
Define
\[
 U(r)=
 \begin{cases}
 v_0(r),&r\ge0,\\
 v_m(r+m),&-m\le r<1-m,\quad m\ge1.
 \end{cases}
\]
The endpoint relations make $U$ continuous.
Every forward restriction of $U$ is obtained by
finitely many time translations and concatenations
of the paths $v_m$.
Proposition~\ref{prop:defect-structure}\textup{(ii), (iv)}
therefore shows that $U$ is a complete defect trajectory.
Strict invariance gives $U(r)\in\cA_*$ for every $r\in\mathbb R$,
so $U$ is bounded and $U(0)=a$.

Conversely, let $U$ be a bounded complete defect trajectory.
The set $B=\{U(t):t\in\mathbb R\}$
is bounded in $H^{-1}_0$.
For every $T>0$,
\[
 U(0)\in\cS_*(T)U(-T)\subseteq\cS_*(T)B.
\]
Applying \eqref{eq:state-attraction} to this bounded
set and using closedness of $\cA_*$ gives
$U(0)\in\cA_*$.
This proves \eqref{eq:defect-kernel-representation}.
This completes the proof.
\end{proof}

\section{Regularity and identification in \texorpdfstring{$L^2$}{L2}}
\label{sec:pullback-l2}

In this section, we establish $L^2$ regularity of bounded
complete defect trajectories and identify the defect
dynamics on $L^2_0$ with the classical evolution.
These results yield $\cA_*=\cA_2$ as subsets of $H^{-1}_0$.

\subsection{\texorpdfstring{$L^2$}{L2} bounds}

We first derive $L^2$ bounds for forward defect trajectories
with $L^2$ initial data and for bounded complete defect
trajectories.
The starting point is the following estimate for the forcing:
\begin{equation}\label{eq:L2-force-weighted-rate}
 \norm{f}_{H^{-1}_\kappa}^2
 =\sum_{j\in\Z}
 \frac{1}{j^2+4\kappa^2}\abs{\widehat f(j)}^2
 \le\frac{1}{4\kappa^2}\norm f_{L^2}^2.
\end{equation}
The forward bounds retain the dependence on the initial
$L^2$ norm, whereas the bounds for complete trajectories
depend only on $\gamma$ and $\norm f_{L^2}$.

\begin{theorem}
\label{thm:pullback-L2}
 Let $f\in L^2_0(\T)$ and $\gamma>0$.
The following $L^2$ regularity properties and bounds hold.
\begin{enumerate}[label=\textup{(\roman*)},leftmargin=2.6em]

\item
Let $u\in\cK_*^+$ with $u(0)\in L^2_0$.
Then $u(t)\in L^2_0$ for every $t\ge0$, and
\begin{equation}\label{eq:forward-L2-bound}
 \sup_{t\ge0}\norm{u(t)}_{L^2}^2
 \le C\left(
 \norm{u(0)}_{L^2}^2
 +\gamma^{-2}\norm f_{L^2}^2
 \right).
\end{equation}

\item
Let $U\in C(\R;H^{-1}_0)$ be a bounded complete
defect trajectory.
Then $U(t)\in L^2_0$ for every $t\in\R$, and
\begin{equation}\label{eq:complete-L2-bound}
 \sup_{t\in\R}\norm{U(t)}_{L^2}^2
 \le C\gamma^{-2}\norm f_{L^2}^2.
\end{equation}
In particular, if $f=0$, then $U\equiv0$.

\end{enumerate}
The constant $C>0$ is absolute.
\end{theorem}

\begin{proof}
\textup{(i)}
By Proposition~\ref{prop:defect-beta}, the determinant
estimates hold along $u$ for all sufficiently large
$\kappa$, uniformly in $t\ge0$.
Taking $s=0$ in \eqref{eq:defect-limit-beta}
and using \eqref{eq:det-weighted-norm}, we have
\begin{align}
 \norm{u(t)}_{H^{-1}_\kappa}^2
 \le C\kappa\alpha_\kappa(u(t)) 
&\le C\norm{u(0)}_{H^{-1}_\kappa}^2
      +C\gamma^{-2}\norm f_{H^{-1}_\kappa}^2\notag\\
 &\le\frac{C}{\kappa^2}
 \left(
 \norm{u(0)}_{L^2}^2
 +\gamma^{-2}\norm f_{L^2}^2
 \right).
 \label{eq:weighted-forward-kappa-minus-two}
\end{align}
The last step uses \eqref{eq:L2-force-weighted-rate}
for both $u(0)$ and $f$.

Since $j^2+4\kappa^2\le5\kappa^2$ for $|j|\le\kappa$, we have
\begin{equation}\label{eq:partial-L2-forward}
 \begin{aligned}
 \sum_{|j|\le\kappa}|\widehat u(t,j)|^2
 \le5\kappa^2\norm{u(t)}_{H^{-1}_\kappa}^2 \le C\left(
 \norm{u(0)}_{L^2}^2
 +\gamma^{-2}\norm f_{L^2}^2
 \right).
 \end{aligned}
\end{equation}
Letting $\kappa\to\infty$ through integers shows that
$u(t)\in L^2_0$.
The bound is independent of $t$, so
\eqref{eq:forward-L2-bound} follows.

\textup{(ii)}
Choose $R>0$ such that
\[
 \sup_{s\in\R}\norm{U(s)}_{H^{-1}}\le R.
\]
Proposition~\ref{prop:defect-beta} gives a common lower
bound $\kappa_R$ for every forward restriction of $U$,
since all initial states have $H^{-1}$ norm at most $R$.
Fix $\kappa\ge\kappa_R$.
For any $s<t$, \eqref{eq:defect-limit-beta} gives
\begin{equation}\label{eq:complete-beta-before-pullback}
 \kappa\alpha_\kappa(U(t))
 \le e^{-c\gamma(t-s)}\kappa\alpha_\kappa(U(s))
     +C\gamma^{-2}\norm f_{H^{-1}_\kappa}^2.
\end{equation}
From \eqref{eq:det-weighted-norm}, we also have
\[
 \sup_{s\in\R}\kappa\alpha_\kappa(U(s))
 \le C\sup_{s\in\R}\norm{U(s)}_{H^{-1}_\kappa}^2
 \le CR^2.
\]
We may therefore let $s\to-\infty$ in \eqref{eq:complete-beta-before-pullback} at this fixed
$\kappa$. Then, by using \eqref{eq:L2-force-weighted-rate}, we obtain
\begin{equation}\label{eq:complete-beta-source-only}
 \kappa\alpha_\kappa(U(t))
 \le C\gamma^{-2}\norm f_{H^{-1}_\kappa}^2
 \le\frac{C}{\gamma^2\kappa^2}\norm f_{L^2}^2.
\end{equation}
As in part \textup{(i)}, the Fourier weights and
\eqref{eq:det-weighted-norm} give
\[
 \begin{aligned}
 \sum_{|j|\le\kappa}|\widehat U(t,j)|^2
 \le5\kappa^2\norm{U(t)}_{H^{-1}_\kappa}^2 \le C\kappa^2
       \bigl[\kappa\alpha_\kappa(U(t))\bigr]\le C\gamma^{-2}\norm f_{L^2}^2.
 \end{aligned}
\]
Letting $\kappa\to\infty$ through integers and taking
the supremum over $t\in\R$ proves
\eqref{eq:complete-L2-bound}.
For $f=0$, this bound gives $U\equiv0$.
This completes the proof.
\end{proof}

\subsection{\texorpdfstring{$L^2$}{L2} continuity and uniqueness}

We now use the scalar Birkhoff equations to establish
strong $L^2$ continuity and identify the trajectories
with the classical evolution.
We begin with the weak continuity implied by the
bounds in Theorem~\ref{thm:pullback-L2}.

Let $I\subset\R$ be an interval with nonempty interior.
Consider a path $u\in C(I;H^{-1}_0)$ satisfying
\eqref{eq:birkhoff-integral} for all $s,t\in I$
with $s\le t$, and assume that
\begin{equation}\label{eq:pathwise-L2-bound-closed-interval}
 \operatorname*{ess\,sup}_{t\in[a,b]}
 \norm{u(t)}_{L^2}<\infty
\end{equation}
for every finite closed subinterval $[a,b]\subseteq I$
with $a<b$.
These assumptions hold for forward defect trajectories
with $L^2$ initial data and for bounded complete defect
trajectories, by Theorems~\ref{thm:pullback-L2}
and~\ref{thm:birkhoff-integral}.

Weak compactness in $L^2$ and continuity in $H^{-1}$
extend the essential bound to every time and give
\begin{equation}\label{eq:weak-L2-time-continuity}
 u\in C_w([a,b];L^2_0),
\end{equation}
including at both endpoints.

Set $z(t)=\Omega(u(t))$.
The $L^2$ bound and \eqref{eq:L2-one-smoothing}
imply that $z$ is bounded in $h^{1/2}$ on $[a,b]$.
By Theorem~\ref{thm:KT}, it is also continuous
in $h^{-1/2}$.
The preceding weak compactness argument therefore gives
$z\in C_w([a,b];h^{1/2})$.

To obtain strong continuity, we remove the linear Airy
phases.
Define the strongly continuous isometric group $E(t)$
on $h^{1/2}$ by
\begin{equation}\label{eq:Airy-phase-h}
 (E(t)z)_n=e^{in^3t}z_n,
 \qquad
 \norm{E(t)z}_{h^{1/2}}=\norm z_{h^{1/2}}.
\end{equation}
For $z\in h^{1/2}$, put $q_z=\Omega^{-1}z\in L^2_0$
and define
\begin{equation}\label{eq:L2-Birkhoff-remainder}
 \mathcal G(z)_n
 =i\omega_n^*(q_z)z_n
  +\bigl(D\Omega(q_z)[f-\gamma q_z]\bigr)_n.
\end{equation}
By \eqref{eq:birkhoff-integral} and
\eqref{eq:frequency-renormalization}, the scalar
coordinates satisfy
\begin{equation}\label{eq:coordinate-split-linear-frequency}
 \dot z_n=in^3z_n+\mathcal G(z)_n
\end{equation}
on the interior of $I$.
Setting $y(t)=E(-t)z(t)$ and integrating, we obtain
\begin{equation}\label{eq:coordinate-interaction-integral}
 y_n(t)-y_n(s)
 =\int_s^t
 \bigl(E(-\tau)\mathcal G(z(\tau))\bigr)_n\,\dd\tau,
 \qquad a\le s\le t\le b.
\end{equation}
The endpoint values are inherited from
\eqref{eq:birkhoff-integral}.

We next show that
$\tau\mapsto E(-\tau)\mathcal G(z(\tau))$
is Bochner integrable with values in $h^{1/2}$.
The local Lipschitz estimates established below will
also give uniqueness for the resulting integral equation.

\begin{lemma}
\label{lem:L2-interaction-Lipschitz}
The map $\mathcal G:h^{1/2}\to h^{1/2}$ defined by
\eqref{eq:L2-Birkhoff-remainder} is locally Lipschitz
and maps bounded sets to bounded sets.
For $t\in\R$, define
\begin{equation}\label{eq:nonautonomous-interaction-field}
 \mathcal G_t(y)=E(-t)\mathcal G(E(t)y).
\end{equation}
Then $\mathcal G_t$ is locally Lipschitz in $y$,
uniformly for $t$ in compact intervals.
Its bounds on bounded sets are uniform for all $t\in\R$.
\end{lemma}

\begin{proof}
By Theorem~\ref{thm:L2-Birkhoff-published}, the map
\[
 z\longmapsto D\Omega(q_z)[f-\gamma q_z]
\]
is analytic and hence locally Lipschitz.
The same theorem shows that
$z\mapsto\omega^*(q_z)$ is locally Lipschitz
from $h^{1/2}$ to $\ell^\infty$.
For $z,w$ near a fixed $z_0\in h^{1/2}$, we therefore have
\begin{align}
 &\norm{
 \bigl(\omega_n^*(q_z)z_n-\omega_n^*(q_w)w_n\bigr)_n
 }_{h^{1/2}}\notag\\
 &\quad\le
 \norm{\omega^*(q_z)}_{\ell^\infty}
 \norm{z-w}_{h^{1/2}}
 +\norm{\omega^*(q_z)-\omega^*(q_w)}_{\ell^\infty}
 \norm w_{h^{1/2}}\notag\\
 &\quad\le C_{z_0}\norm{z-w}_{h^{1/2}}.
 \label{eq:frequency-product-Lipschitz}
\end{align}
Combining the two terms proves local Lipschitz continuity
of $\mathcal G$.
The local constants may depend on $f$ and $\gamma$.

Now fix $r>0$.
Since $\Omega^{-1}$ maps bounded sets to bounded sets,
there is $R>0$ such that
$\norm{q_z}_{L^2}\le R$ whenever
$\norm z_{h^{1/2}}\le r$.
Using \eqref{eq:L2-frequency-tail} and
Lemma~\ref{lem:Birkhoff-derivative-balls}, we get
\[
 \begin{aligned}
 \norm{\mathcal G(z)}_{h^{1/2}}
 &\le
 \norm{\omega^*(q_z)}_{\ell^\infty}\norm z_{h^{1/2}}
 +\norm{D\Omega(q_z)[f-\gamma q_z]}_{h^{1/2}}\\
 &\le C_R\bigl(r+\norm f_{L^2}+\gamma R\bigr).
 \end{aligned}
\]
Thus $\mathcal G$ is bounded on bounded sets.

To prove the local uniformity in time, fix
$y_0\in h^{1/2}$ and a compact interval $J\subset\R$.
Strong continuity of $E(t)$ makes the orbit
$\{E(t)y_0:t\in J\}$ compact.
A finite-cover argument using local Lipschitz continuity
of $\mathcal G$ and the isometries $E(t)$ gives
$\delta>0$ and $C_{J,y_0}<\infty$ such that
\[
 \norm{\mathcal G_t(y)-\mathcal G_t(w)}_{h^{1/2}}
 \le C_{J,y_0}\norm{y-w}_{h^{1/2}}
\]
whenever $t\in J$ and $y,w$ lie in the $\delta$-ball
about $y_0$.

Finally, $E(t)$ preserves every centered ball, so
\[
 \sup_{t\in\R}\sup_{\norm y_{h^{1/2}}\le r}
 \norm{\mathcal G_t(y)}_{h^{1/2}}
 =
 \sup_{\norm z_{h^{1/2}}\le r}
 \norm{\mathcal G(z)}_{h^{1/2}}
 <\infty
 \qquad(r>0).
\]
This completes the proof.
\end{proof}

We now establish strong $L^2$ continuity and identify
the trajectories with the classical evolution.
The conclusion includes any endpoints belonging to
the time interval.

\begin{theorem}
\label{thm:L2-coordinate-identification}
Let $I\subset\R$ be an interval with nonempty interior,
and let $u\in C(I;H^{-1}_0)$.
Assume \eqref{eq:pathwise-L2-bound-closed-interval}
on every finite closed subinterval $[a,b]\subseteq I$
with $a<b$, and suppose that
\eqref{eq:birkhoff-integral} holds for every $n\ge1$
and every $s,t\in I$ with $s\le t$.
Then
\begin{equation}\label{eq:L2-strong-continuity-conclusion}
 u\in C(I;L^2_0),
\end{equation}
and
\[
 u(t)=S_2(t-s)u(s),
 \qquad s,t\in I,\quad s\le t.
\]
In particular, for $I=[0,\infty)$,
\begin{equation}\label{eq:L2-initial-trace}
 \lim_{t\downarrow0}\norm{u(t)-u(0)}_{L^2}=0,
 \qquad
 u(t)=S_2(t)u(0)\quad(t\ge0).
\end{equation}
Moreover, $u$ is a standard distributional solution:
for every $\varphi\in C_c^\infty(I^\circ\times\T)$,
\begin{equation}\label{eq:standard-distributional-formulation}
 \int_{I^\circ\times\T}
 \left[
 -u\,\partial_t\varphi-u\,\partial_x^3\varphi
 +3u^2\partial_x\varphi+\gamma u\varphi-f\varphi
 \right]\dd\mu(x)\dd t=0.
\end{equation}
\end{theorem}

\begin{proof}
We first establish strong continuity in $L^2$.
Fix $[a,b]\subseteq I$ with $a<b$, and put
$z=\Omega(u)$ and $y(t)=E(-t)z(t)$.
From the preceding discussion, $z$ is bounded and
weakly continuous in $h^{1/2}$, while $y$ satisfies
\eqref{eq:coordinate-interaction-integral}.
Set
\[
 H(\tau)=E(-\tau)\mathcal G(z(\tau))
        =\mathcal G_\tau(y(\tau)).
\]
By Lemma~\ref{lem:L2-interaction-Lipschitz}, we have
\[
 \sup_{\tau\in[a,b]}\norm{H(\tau)}_{h^{1/2}}\le C.
\]
Each coordinate of $H$ is continuous by the
$H^{-1}$ continuity of $u$ and Theorem~\ref{thm:KT}.
Its finite truncations are therefore continuous and
converge pointwise in $h^{1/2}$ to $H$.
It follows that $H$ is strongly measurable and,
being bounded, belongs to $L^1([a,b];h^{1/2})$.

Applying the coordinate functionals to the Bochner
integral and using
\eqref{eq:coordinate-interaction-integral},
one obtains
\begin{equation}\label{eq:vector-interaction-integral}
 y(t)-y(s)
 =\int_s^t\mathcal G_\tau(y(\tau))\,\dd\tau
 \quad\text{in }h^{1/2},
 \qquad a\le s\le t\le b.
\end{equation}
In particular, we have
\begin{equation}\label{eq:interaction-closed-interval-modulus}
 \norm{y(t)-y(s)}_{h^{1/2}}
 \le C|t-s|,
 \qquad s,t\in[a,b].
\end{equation}
Thus $y$ is Lipschitz on $[a,b]$.
Since $E(t)$ is strongly continuous, it is easy to see
that $z(t)=E(t)y(t)$ is continuous in $h^{1/2}$.
The continuity of $\Omega^{-1}:h^{1/2}\to L^2_0$
then gives $u\in C([a,b];L^2_0)$, including at both
endpoints.
As $[a,b]$ was arbitrary, this proves
\eqref{eq:L2-strong-continuity-conclusion}.

It remains to identify $u$ with the classical evolution.
Fix $[s,b]\subseteq I$ with $s<b$, and define
\[
 v(t)=S_2(t-s)u(s),\qquad t\ge s.
\]
The smooth approximation property in
Theorem~\ref{thm:classical-L2} shows that
$r\mapsto v(s+r)$ is a primitive trajectory.
We therefore have the scalar coordinate identities
for $v$ by Proposition~\ref{prop:primitive-inherited}.
Applying the first part of the proof to $v$, we find
that $\widetilde y(t)=E(-t)\Omega(v(t))$
satisfies \eqref{eq:vector-interaction-integral}
on $[s,b]$, with $\widetilde y(s)=y(s)$.

By continuity of the two paths and
Lemma~\ref{lem:L2-interaction-Lipschitz}, we may choose
$0<h\le b-s$ so that both remain in a neighborhood
on which $\mathcal G_\tau$ is Lipschitz with a common
constant $L$, uniformly for $\tau\in[s,s+h]$.
Subtracting their integral equations, we obtain
\[
 \norm{y(r)-\widetilde y(r)}_{h^{1/2}}
 \le L\int_s^r
 \norm{y(\tau)-\widetilde y(\tau)}_{h^{1/2}}\,\dd\tau,
 \qquad s\le r\le s+h.
\]
By Gronwall's inequality, $y(r)=\widetilde y(r)$
for every $r\in[s,s+h]$.
The same local argument, together with continuity,
extends this identity to $[s,b]$.
Inverting the coordinate transform, we get
\[
 u(t)=v(t)=S_2(t-s)u(s),
 \qquad s\le t\le b.
\]
Since $[s,b]$ was arbitrary, the evolution identity
holds throughout $I$.
For $I=[0,\infty)$, taking $s=0$ and using strong
continuity at zero gives \eqref{eq:L2-initial-trace}.

Finally, let $\varphi\in C_c^\infty(I^\circ\times\T)$
and choose $s\in I$ earlier than its time support.
On this support, $u$ agrees with the classical solution
$S_2(t-s)u(s)$.
The weak formulation for that solution gives
\eqref{eq:standard-distributional-formulation}.
This completes the proof.
\end{proof}

We can now identify the trajectories in
Theorem~\ref{thm:pullback-L2} with the classical evolution
and recover the usual quadratic nonlinearity.

\begin{corollary}
\label{cor:complete-standard-distributional}
Let $u\in\cK_*^+$ satisfy $u(0)\in L^2_0$, and let
$U$ be a bounded complete defect trajectory.
Then
\[
 u\in C([0,\infty);L^2_0),
 \qquad
 u(t)=S_2(t-s)u(s)\quad(0\le s\le t).
\]
Moreover, $U$ is a bounded complete trajectory of $S_2$:
\[
 U\in C(\R;L^2_0),
 \qquad
 U(t)=S_2(t-s)U(s)\quad(s\le t),
\]
with $\sup_{t\in\R}\norm{U(t)}_{L^2}<\infty$.
The nonlinear distributions satisfy
\begin{equation}\label{eq:closed-product-equals-standard-complete}
 \begin{aligned}
 \mathfrak N[u]&=\partial_x(u^2)
 &&\text{in }\mathcal D'((0,\infty)\times\T),\\
 \mathfrak N[U]&=\partial_x(U^2)
 &&\text{in }\mathcal D'(\R\times\T).
 \end{aligned}
\end{equation}
\end{corollary}

\begin{proof}
From Theorem~\ref{thm:pullback-L2}, we have
\[
 \sup_{t\ge0}\norm{u(t)}_{L^2}<\infty,
 \qquad
 \sup_{t\in\R}\norm{U(t)}_{L^2}<\infty.
\]
The scalar coordinate identities for $u$ follow from
Theorem~\ref{thm:birkhoff-integral}.
For $U$, completeness allows us to apply the same
theorem to every forward restriction.
Thus Theorem~\ref{thm:L2-coordinate-identification}
applies on $[0,\infty)$ and $\R$, respectively,
giving the asserted strong continuity and semigroup
identities.
Together with the second bound above, this makes $U$
a bounded complete trajectory of $S_2$.

It remains to identify the nonlinear distributions.
The $L^2$ continuity ensures that $u^2$ and $U^2$
are locally integrable on their respective time domains.
For $U$, apply \eqref{eq:pathwise-defect-nonlinearity}
to its forward restrictions.
The local expressions agree on overlaps, since each
is given by the same linear expression in $U$.
They therefore define $\mathfrak N[U]$ on $\R\times\T$.

Using \eqref{eq:pathwise-defect-nonlinearity} and
\eqref{eq:standard-distributional-formulation},
we obtain
\[
 \begin{aligned}
 \mathfrak N[u]
 &=\frac13\bigl(
 \partial_tu+\partial_x^3u+\gamma u-f
 \bigr)
 =\partial_x(u^2),\\
 \mathfrak N[U]
 &=\frac13\bigl(
 \partial_tU+\partial_x^3U+\gamma U-f
 \bigr)
 =\partial_x(U^2).
 \end{aligned}
\]
The first identity holds in
$\mathcal D'((0,\infty)\times\T)$ and the second
in $\mathcal D'(\R\times\T)$.
This completes the proof.
\end{proof}
\subsection{Identification of the attractors}

We now complete the proofs of the two main theorems on defect
dynamics. After collecting the results for the evolution on
$H^{-1}_0$, we identify its $L^2$ trajectories and global attractor
with their classical counterparts.

\begin{proof}[Proof of Theorem~\ref{thm:intro-defect}]
By Theorem~\ref{thm:strict-multiflow},
$\{\cS_*(t)\}_{t\ge0}$ is a strict multivalued semiflow
on $H^{-1}_0$ with nonempty values, and
$\Graph\cS_*(t)$ is closed for every fixed $t\ge0$.
Moreover, Theorem~\ref{thm:multivalued-attractor}
provides a nonempty compact global attractor $\cA_*$
with the representation
\eqref{eq:defect-kernel-representation}.
Uniqueness follows from strict invariance and
attraction of bounded sets.
This completes the proof.
\end{proof}

\begin{proof}[Proof of Theorem~\ref{thm:intro-L2-identification}]
We first identify the forward trajectories.
Let $u_0\in L^2_0$.
By Proposition~\ref{prop:primitive-existence},
there is a primitive trajectory starting at $u_0$.
Since $\cK_0^+\subset\cK_*^+$,
Corollary~\ref{cor:complete-standard-distributional}
identifies this path with $S_2(\cdot)u_0$.
The same corollary applies to every defect trajectory
with initial state $u_0$. Hence
\[
 S_2(\cdot)u_0\in\cK_0^+,
 \qquad
 \{u\in\cK_*^+:u(0)=u_0\}
 =\{S_2(\cdot)u_0\},
\]
which proves
\eqref{eq:defect-forward-L2-identification}.

Next let us consider bounded complete trajectories.
Let $U$ be a bounded complete defect trajectory.
From Theorem~\ref{thm:pullback-L2}\textup{(ii)}
and Corollary~\ref{cor:complete-standard-distributional},
one has
\[
 \begin{gathered}
 \sup_{t\in\R}\norm{U(t)}_{L^2}^2
 \le C\gamma^{-2}\norm f_{L^2}^2,\\
 U\in C(\R;L^2_0),
 \qquad
 U(t)=S_2(t-s)U(s)\quad(s\le t).
 \end{gathered}
\]
The first bound proves \eqref{eq:intro-L2},
and these properties show that $U$ is a bounded
complete trajectory of $S_2$.

Conversely, let $V$ be a bounded complete trajectory
of $S_2$.
The continuous embedding $L^2_0\hookrightarrow H^{-1}_0$
ensures that $V$ is continuous and bounded in $H^{-1}_0$.
For every $s\in\R$, the forward identification above gives
\[
 (\mathcal T_sV)|_{[0,\infty)}
 =S_2(\cdot)V(s)\in\cK_0^+\subset\cK_*^+.
\]
Thus $V$ is a bounded complete defect trajectory.
This establishes the identification in both directions.

Finally, the kernel representations
\eqref{eq:defect-kernel-representation} and
\eqref{eq:L2-attractor-kernel} give
\[
 \begin{aligned}
 \cA_*
 &=\{U(0):
 U\text{ is a bounded complete defect trajectory}\}\\
 &=\{V(0):
 V\text{ is a bounded complete trajectory of }S_2\}.
 \end{aligned}
\]
Consequently,
\begin{equation}\label{eq:spectral-attractor-classical}
 \cA_*=\cA_2
 \qquad\text{as subsets of }H^{-1}_0.
\end{equation}

The same argument also identifies the primitive trajectory attractor.
Indeed, every bounded complete primitive trajectory is a bounded
complete defect trajectory, while every forward restriction of a
bounded complete trajectory of $S_2$ belongs to $\cK_0^+$.
The representation in Theorem~\ref{thm:intro-primitive}, together
with \eqref{eq:L2-attractor-kernel}, therefore gives
\[
 \cA=\{S_2(\cdot)a:a\in\cA_2\},
 \qquad
 \{u(0):u\in\cA\}=\cA_*=\cA_2.
\]

When $f=0$, \eqref{eq:complete-L2-bound}
implies that every bounded complete defect trajectory
is identically zero.
The kernel representation and nonemptiness of $\cA_*$
then yield $\cA_*=\cA_2=\{0\}$.
This completes the proof.
\end{proof}

\section*{Declaration of AI use}

During the preparation of this manuscript, the author
used ChatGPT and Codex (OpenAI) to assist with mathematical
exploration, proof development, manuscript organization,
and English-language revision.
The AI-assisted material was reviewed and revised by
the author, who takes full responsibility for the
mathematical content, references, and final manuscript.

\section*{Acknowledgments}

  M.S is partially supported by the National Natural Science Foundation of China, grant numbers: 12101629 and 12371123.

\bigskip
\noindent\begin{minipage}{\textwidth}
\footnotesize
\textsc{College of Science, Minzu University of China,
Beijing 100081, P. R. China}\par
\textit{E-mail address:}
\href{mailto:smj@muc.edu.cn}{\texttt{smj@muc.edu.cn}}
\end{minipage}


\begin{thebibliography}{99}
\setlength{\itemsep}{2pt}

\bibitem{BabinIlyinTiti2011}
A. V. Babin, A. A. Ilyin, and E. S. Titi,
\emph{On the regularization mechanism for the periodic Korteweg--de Vries equation},
Comm. Pure Appl. Math. \textbf{64} (2011), 591--648.

\bibitem{Ball1997}

J. M. Ball, \emph{Continuity properties and global attractors of generalized
semiflows and the Navier--Stokes equations},
J. Nonlinear Sci. \textbf{7} (1997), 475--502.


\bibitem{BonaZhang1996}
J. L. Bona and B.-Y. Zhang, \emph{The initial-value problem for the forced
Korteweg--de Vries equation}, Proc. Roy. Soc. Edinburgh Sect. A
\textbf{126} (1996), 571--598.

\bibitem{Bourgain1993KdV}
J. Bourgain, \emph{Fourier transform restriction phenomena for certain lattice
subsets and applications to nonlinear evolution equations. II. The KdV-equation},
Geom. Funct. Anal. \textbf{3} (1993), 209--262.

\bibitem{Bourgain1997}

J. Bourgain, \emph{Periodic Korteweg de Vries equation with measures as initial data},
Selecta Math. (N.S.) \textbf{3} (1997), 115--159.


\bibitem{ChepyzhovVishik1997}

V. V. Chepyzhov and M. I. Vishik, \emph{Evolution equations and their trajectory
attractors}, J. Math. Pures Appl. \textbf{76} (1997), 913--964.


\bibitem{ChepyzhovVishik2002}

V. V. Chepyzhov and M. I. Vishik, \emph{Attractors for Equations of Mathematical
Physics}, Amer. Math. Soc. Colloq. Publ. 49, AMS, Providence, RI, 2002.


\bibitem{ChepyzhovVishikZelik2011}
V. V. Chepyzhov, M. I. Vishik, and S. V. Zelik, \emph{Strong trajectory attractors
for dissipative Euler equations}, J. Math. Pures Appl. \textbf{96} (2011), 395--407.

\bibitem{Cheskidov2009}

A. Cheskidov, \emph{Global attractors of evolutionary systems},
J. Dynam. Differential Equations \textbf{21} (2009), 249--268.


\bibitem{CKSTT2003}
J. Colliander, M. Keel, G. Staffilani, H. Takaoka, and T. Tao,
\emph{Sharp global well-posedness for KdV and modified KdV on $\mathbb R$ and
$\mathbb T$}, J. Amer. Math. Soc. \textbf{16} (2003), 705--749.

\bibitem{CKSTT2004}
J. Colliander, M. Keel, G. Staffilani, H. Takaoka, and T. Tao,
\emph{Multilinear estimates for periodic KdV equations, and applications},
J. Funct. Anal. \textbf{211} (2004), 173--218.

\bibitem{ErdoganTzirakis2012}
M. B. Erdo\u{g}an and N. Tzirakis, \emph{Global smoothing for the periodic KdV
evolution}, Int. Math. Res. Not. \textbf{2013} (2013), 4589--4614.

\bibitem{ErdoganTzirakis2013}
M. B. Erdo\u{g}an and N. Tzirakis, \emph{Long time dynamics for forced and weakly
damped KdV on the torus}, Commun. Pure Appl. Anal. \textbf{12} (2013), 2669--2684.

\bibitem{Gardner1967}

C. S. Gardner, J. M. Greene, M. D. Kruskal, and R. M. Miura,
\emph{Method for solving the Korteweg--deVries equation},
Phys. Rev. Lett. \textbf{19} (1967), 1095--1097.


\bibitem{Ghidaglia1988}
J.-M. Ghidaglia, \emph{Weakly damped forced Korteweg--de Vries equations behave
as a finite dimensional dynamical system in the long time},
J. Differential Equations \textbf{74} (1988), 369--390.

\bibitem{Ghidaglia1994}
J.-M. Ghidaglia, \emph{A note on the strong convergence towards attractors of
damped forced KdV equations}, J. Differential Equations \textbf{110} (1994), 356--359.

\bibitem{Goubet2000}
O. Goubet, \emph{Asymptotic smoothing effect for weakly damped forced
Korteweg--de Vries equations}, Discrete Contin. Dyn. Syst. \textbf{6} (2000), 625--644.

\bibitem{Goubet2018}
O. Goubet, \emph{Analyticity of the global attractor for damped forced periodic
Korteweg--de Vries equation}, J. Differential Equations \textbf{264} (2018), 3052--3066.

\bibitem{GoubetMolinet2009}

O. Goubet and L. Molinet,
\emph{Global attractor for weakly damped nonlinear Schr\"odinger equations
in $L^2(\mathbb R)$}, Nonlinear Anal. \textbf{71} (2009), nos.~1--2, 317--320.
\href{https://doi.org/10.1016/j.na.2008.10.078}{doi:10.1016/j.na.2008.10.078}.


\bibitem{GoubetRosa2002}

O. Goubet and R. Rosa, \emph{Asymptotic smoothing and the global attractor of
a weakly damped KdV equation on the real line},
J. Differential Equations \textbf{185} (2002), 25--53.


\bibitem{Goyal2021}
P. Goyal, \emph{Global attractor for weakly damped, forced mKdV equation below
energy space}, Nagoya Math. J. \textbf{241} (2021), 171--203.

\bibitem{Guo2009}

Z. Guo, \emph{Global well-posedness of Korteweg--de Vries equation in
$H^{-3/4}(\mathbb R)$}, J. Math. Pures Appl. \textbf{91} (2009), 583--597.


\bibitem{Hale1988}

J. K. Hale, \emph{Asymptotic Behavior of Dissipative Systems},
Math. Surveys Monogr. 25, AMS, Providence, RI, 1988.


\bibitem{HanKloeden2023}
X. Han and P. Kloeden, \emph{Dissipative Lattice Dynamical Systems},
Interdiscip. Math. Sci. 22, World Scientific, Singapore, 2023.
\href{https://doi.org/10.1142/13182}{doi:10.1142/13182}.

\bibitem{KappelerMoehrTopalov2005}
T. Kappeler, C. M\"ohr, and P. Topalov, \emph{Birkhoff coordinates for KdV on
phase spaces of distributions}, Selecta Math. (N.S.) \textbf{11} (2005), 37--98.

\bibitem{KappelerMolnar2018}
T. Kappeler and J.-C. Molnar, \emph{On the wellposedness of the KdV/KdV2
equations and their frequency maps},
Ann. Inst. H. Poincar\'e C Anal. Non Lin\'eaire \textbf{35} (2018), 101--160.

\bibitem{KappelerPoschel2003}

T. Kappeler and J. P\"oschel, \emph{KdV \& KAM},
Ergebnisse der Mathematik und ihrer Grenzgebiete 45, Springer, Berlin, 2003.


\bibitem{KappelerSchaadTopalov2013}
T. Kappeler, B. Schaad, and P. Topalov, \emph{Qualitative features of periodic
solutions of KdV}, Comm. Partial Differential Equations \textbf{38} (2013), 1626--1673.

\bibitem{KappelerSerierTopalov2008}

T. Kappeler, F. Serier, and P. Topalov, \emph{On the symplectic phase space of
KdV}, Proc. Amer. Math. Soc. \textbf{136} (2008), 1691--1698.


\bibitem{KappelerTopalovRiccati2005}

T. Kappeler and P. Topalov, \emph{Riccati map on $L^2_0(\mathbb T)$ and its
applications}, J. Math. Anal. Appl. \textbf{309} (2005), 544--566.


\bibitem{KappelerTopalov2006}
T. Kappeler and P. Topalov, \emph{Global wellposedness of KdV in
$H^{-1}(\mathbb T,\mathbb R)$}, Duke Math. J. \textbf{135} (2006), 327--360.

\bibitem{Kato1979}
T. Kato, \emph{On the Korteweg--de Vries equation},
Manuscripta Math. \textbf{28} (1979), 89--99.

\bibitem{KenigPonceVega1996}
C. E. Kenig, G. Ponce, and L. Vega, \emph{A bilinear estimate with applications
to the KdV equation}, J. Amer. Math. Soc. \textbf{9} (1996), 573--603.

\bibitem{KillipVisan2019}
R. Killip and M. Vi\c{s}an, \emph{KdV is well-posed in $H^{-1}$},
Ann. of Math. \textbf{190} (2019), 249--305.

\bibitem{KillipVisanZhang2018}
R. Killip, M. Vi\c{s}an, and X. Zhang, \emph{Low regularity conservation laws
for integrable PDE}, Geom. Funct. Anal. \textbf{28} (2018), 1062--1090.
\href{https://arxiv.org/abs/1708.05362}{arXiv:1708.05362}.

\bibitem{Kishimoto2009}

N. Kishimoto, \emph{Well-posedness of the Cauchy problem for the Korteweg--de
Vries equation at the critical regularity},
Differential Integral Equations \textbf{22} (2009), 447--464.


\bibitem{Lax1968}

P. D. Lax, \emph{Integrals of nonlinear equations of evolution and solitary
waves}, Comm. Pure Appl. Math. \textbf{21} (1968), 467--490.


\bibitem{LiWu2009}
Y. Li and Y. Wu, \emph{Global attractor for weakly damped forced KdV equation
in low regularity on $\mathbb T$}, preprint, 2009,
\href{https://arxiv.org/abs/0910.4652}{arXiv:0910.4652}.

\bibitem{McConnell2023}
R. McConnell, \emph{Global attractor for the periodic generalized
Korteweg--de Vries equation through smoothing},
Discrete Contin. Dyn. Syst. Ser. B \textbf{28} (2023), 1133--1158.

\bibitem{McKeanTrubowitz1976}

H. P. McKean and E. Trubowitz, \emph{Hill's operator and hyperelliptic
function theory in the presence of infinitely many branch points},
Comm. Pure Appl. Math. \textbf{29} (1976), 143--226.


\bibitem{MelnikValero1998}
V. S. Melnik and J. Valero, \emph{On attractors of multivalued semi-flows and
differential inclusions}, Set-Valued Anal. \textbf{6} (1998), 83--111.

\bibitem{MelnikValero2008}

V. S. Melnik and J. Valero, \emph{Addendum to ``On attractors of multivalued
semiflows and differential inclusions''}, Set-Valued Anal. \textbf{16} (2008), 507--509.


\bibitem{Molinet2009Attractor}
L. Molinet, \emph{Global attractor and asymptotic smoothing effects for
the weakly damped cubic Schr\"odinger equation in $L^2(\mathbb T)$},
Dyn. Partial Differ. Equ. \textbf{6} (2009), 15--34.
\href{https://doi.org/10.4310/DPDE.2009.v6.n1.a2}{doi:10.4310/DPDE.2009.v6.n1.a2}.

\bibitem{Molinet2012}
L. Molinet, \emph{Sharp ill-posedness results for the KdV and mKdV equations
on the torus}, Adv. Math. \textbf{230} (2012), 1895--1930.

\bibitem{MolinetVento2011}

L. Molinet and S. Vento, \emph{Sharp ill-posedness and well-posedness results
for the KdV--Burgers equation: the real line case},
Ann. Sc. Norm. Super. Pisa Cl. Sci. \textbf{10} (2011), 531--560.


\bibitem{Robinson2001}

J. C. Robinson, \emph{Infinite-Dimensional Dynamical Systems},
Cambridge University Press, Cambridge, 2001.


\bibitem{Simon2005}
B. Simon, \emph{Trace Ideals and Their Applications}, second ed.,
Math. Surveys Monogr. 120, AMS, Providence, RI, 2005.

\bibitem{Temam1997}

R. Temam, \emph{Infinite-Dimensional Dynamical Systems in Mechanics and
Physics}, second ed., Springer, New York, 1997.


\bibitem{Tsugawa2004}
K. Tsugawa, \emph{Existence of the global attractor for weakly damped, forced
KdV equation on Sobolev spaces of negative index},
Commun. Pure Appl. Anal. \textbf{3} (2004), 301--318.

\bibitem{VishikChepyzhov2011}
M. I. Vishik and V. V. Chepyzhov, \emph{Trajectory attractors of equations of
mathematical physics}, Russian Math. Surveys \textbf{66} (2011), 637--731.

\bibitem{Wang2014BBM}

M. Wang, \emph{Long time dynamics for a damped Benjamin--Bona--Mahony
equation in low regularity spaces}, Nonlinear Anal. \textbf{105}
(2014), 134--144.
\href{https://doi.org/10.1016/j.na.2014.04.013}{doi:10.1016/j.na.2014.04.013}.


\bibitem{Wang2015gKdV}
M. Wang, \emph{Global attractor for weakly damped gKdV equations in higher
Sobolev spaces}, Discrete Contin. Dyn. Syst. \textbf{35} (2015), 3799--3825.
\href{https://doi.org/10.3934/dcds.2015.35.3799}{doi:10.3934/dcds.2015.35.3799}.

\bibitem{Wang2024BBM}
M. Wang, \emph{Global attractor for the damped BBM equation in the sharp
low regularity space}, Z. Angew. Math. Phys. \textbf{75} (2024), Art.~142.
\href{https://doi.org/10.1007/s00033-024-02288-7}{doi:10.1007/s00033-024-02288-7}.

\bibitem{Wang2024BBM3D}
M. Wang, \emph{Well Posedness and Global Attractors for the 3D Periodic
BBM Equation Below the Energy Space},
J. Dynam. Differential Equations \textbf{36} (2024), no.~4, 3599--3621.
\href{https://doi.org/10.1007/s10884-023-10253-7}{doi:10.1007/s10884-023-10253-7}.

\bibitem{WangLiu2018BBM}
M. Wang and A. Liu, \emph{Dynamics of the BBM equation with a distribution
force in low regularity spaces}, Topol. Methods Nonlinear Anal.
\textbf{51} (2018), 91--109.
\href{https://doi.org/10.12775/TMNA.2017.058}{doi:10.12775/TMNA.2017.058}.

\bibitem{Yang2011}
X. Yang, \emph{Global attractor for the weakly damped forced KdV equation
in Sobolev spaces of low regularity},
NoDEA Nonlinear Differential Equations Appl. \textbf{18} (2011), 273--285.

\bibitem{Zelik2023}

S. V. Zelik, \emph{Attractors. Then and now},
Russian Math. Surveys \textbf{78} (2023), 635--777.


\end{thebibliography}
\end{document}